\documentclass[12pt]{article}

\usepackage{amsfonts}

\usepackage{amsmath,amssymb,amstext,amsthm,amscd,mathtools,graphicx,float,caption,subcaption,todonotes,tikz,tikz-cd,hyperref}

\usetikzlibrary{matrix,arrows,decorations.pathmorphing}

\theoremstyle{plain}
\newtheorem{thm}{Theorem}[section] 
\newtheorem{pro}[thm]{Proposition}
\newtheorem{lem}[thm]{Lemma}

\newtheorem{con}{Conjecture}

\theoremstyle{definition}
\newtheorem*{defi}{Definition}
\newtheorem{remark}[thm]{Remark}

\theoremstyle{remark}

\def\ds{\displaystyle}
\def\pf{{\em Proof.}\ \,}

\def\su{\subseteq}

\def\({\left(}
\def\){\right)}

\def\fun{\rightarrow}
\def\lfun{\longrightarrow}
\def\sur{\twoheadrightarrow}

\def\bC{{\mathbb{C}}}

\def\bR{{\mathbb{R}}}

\def\cA{{\mathcal A}}

\def\cF{{\mathcal F}}

\def\cI{{\mathcal I}}

\def\cO{{\mathcal O}}

\def\cS{{\mathcal S}}
\def\cT{{\mathcal T}}

\def\cW{{\mathcal W}}
\def\cZ{{\mathcal Z}}

\def\ch{{\textrm{ch}}}

\def\Coh{{\textrm{Coh\,}}}

\def\Hom{{\textrm{Hom\,}}}

\def\Im{{\textrm{Im\,}}}

\def\Real{{\textrm{Re\,}}}
\def\rk{{\textrm{rk}}}

\newcommand{\stabdivS}{\text{Stab}_{div}(S)}
\newcommand{\DS}{D(S)}
\newcommand{\Kl}{\underline{K}}

\newcommand{\Jl}{\overline{J}}

\begin{document}

\title{Bridgeland Stability of Line Bundles on Surfaces}

\author{Daniele Arcara and Eric Miles}

\date{}

\maketitle

\begin{abstract}
We study the Bridgeland stability of line bundles on surfaces with respect to certain Bridgeland stability conditions determined by divisors.
Given a smooth projective surface $S$, we show that a line bundle $L$ is always Bridgeland stable for those stability conditions if there are no curves $C\su S$ of negative self-intersection.
When a curve $C$ of negative self-intersection is present, $L$ is destabilized by $L(-C)$ for some stability conditions.
We conjecture that line bundles of the form $L(-C)$ are the only objects that can destabilize $L$ and that torsion sheaves of the form $L(C)|_C$ are the only objects that can destabilize $L[1]$.
We prove our conjecture in several cases, in particular for Hirzebruch surfaces.
\end{abstract} 

\tableofcontents



\section{Introduction}

Let $S$ be a smooth projective surface. In this paper, we study Bridgeland stability for line bundles on $S$ using the geometric Bridgeland stability conditions introduced in \cite{ABL} (see Section \ref{stabcond} for a precise definition). They are stability conditions for complexes of sheaves in the bounded derived category $D^b(\Coh S)$. Line bundles are always Mumford slope-stable because their only subobjects are ideal sheaves; however, the situation is less constrained in the derived setting. For example, in the Abelian subcategories that we consider, a subobject of a line bundle is a sheaf but may a priori have arbitrarily high rank. The quotient is possibly a two-term complex. 

One might still expect line bundles to always be Bridgeland stable, and this is correct if $S$ has no curves $C$ of negative self-intersection (see the first part of Theorem \ref{mainthm} below).
However, if there exists a curve $C$ on $S$ of negative self-intersection, then there exist Bridgeland stability conditions for which $L(-C)$ destabilizes $L$, as well as Bridgeland stability conditions for which $L(C)|_C$ destabilizes $L[1]$.
We make the following conjecture.
\begin{con}\label{mainconj}
Given a surface $S$ and a stability condition $\sigma_{H,D}$ as in \cite{ABL},
\begin{itemize}
\item
the only objects that could destabilize a line bundle $L$ are line bundles of the form $L(-C)$ for a curve $C$ of negative self-intersection, and
\item
the only objects that could destabilize $L[1]$ are torsion sheaves of the form $L(C)|_C$ for a curve $C$ of negative self-intersection.
\end{itemize}
\end{con}

The goal of this paper is to prove the conjecture in several cases and provide evidence for the conjecture in others. Specifically, we prove the following.

\begin{thm}\label{mainthm}
The conjecture is true in the following cases: 
\begin{itemize}
\item
If $S$ does not have any curves of negative self-intersection.
\item
If the Picard rank of $S$ is $2$ and there exists only one irreducible curve of negative self-intersection.
\end{itemize}
\end{thm}

In particular, the conjecture is true for Hirzebruch surfaces because they have Picard rank 2 and only one irreducible curve of negative self-intersection (see, e.g., \cite{Beauville}).

We cite Propositions \ref{LH} and \ref{conjevidencePicRk2} as further evidence for our conjecture in general. Proposition \ref{LH} establishes some structure for the walls of an actually destabilizing subobject of a line bundle on surfaces of any Picard rank, and Proposition \ref{conjevidencePicRk2} proves a stronger version of the conjecture for a subset of stability conditions when $S$ has Picard rank 2 and two irreducible curves of negative self-intersection.

In \cite{ABL}, the stability of line bundles is proven for stability conditions $\sigma_{D,H}$ with $D=sH$ and is utilized in \cite{ABCH} and \cite{biratgeomhilbsurfs} when classifying destabilizing walls for ideal sheaves of points on surfaces. When $D\neq sH$, the more algebraic proof of \cite{ABL} using the Bogomolov inequality and Hodge Index Theorem fails, and new techniques are required. We use Theorem 3.1 (Bertram's Nested Wall Theorem) from \cite{MaciociaWalls} and Lemma 6.3 (which we refer to as Bertram's Lemma) from \cite{ABCH} along with an analysis of the relative geometry of relevant walls in certain three-dimensional slices of the space of stability conditions.

We begin in Section \ref{BSCs} by introducing Bridgeland stability conditions and the stability conditions $\sigma_{D,H}$ of interest, as well as important slices of the space of stability conditions. In Section \ref{reduction}, we present an action by line bundles that allows us to consider only $\cO_S$ in our question concerning the stability of line bundles. In Section \ref{prelims}, we consider the basic structure of subobjects of $\cO_S$ and present two previously obtained results that will serve as important tools in the remainder of this paper. In Section \ref{stabofO}, we prove our main results, but we first consider the rank 1 subobjects of $\cO_S$. The rank 1 subobjects form the base case for our main results, all of which use induction. The case of $\cO_S[1]$ is then primarily completed using duality, which allows us to use our results for $\cO_S$ except when $D.H=0$.


\section*{Acknowledgments}
We would like to thank Arend Bayer, Cristian Martinez, and especially Aaron Bertram and Renzo Cavalieri for many useful discussions.
We would also like to thank the referee for many useful comments and corrections.


\section{Bridgeland stability conditions}
\label{BSCs}

Bridgeland stability conditions (introduced in \cite{stabcondsontricats}) provide a concept of stability on the derived category of a variety.
As with Mumford-slope stability, we may deform our stability conditions (Bridgeland showed that the space of all stability conditions is a complex manifold), and the stability of objects can change.
We first introduce these stability conditions in general; then, we restrict our attention to surfaces in the next section.
Our goal is to study the Bridgeland stability of line bundles.


\subsection{General definition of Bridgeland stability conditions}

Let $X$ be a smooth projective variety, $D(X)$ the bounded derived category of coherent sheaves on $X$, $K(X)$ its Grothendieck group, and $K_{num}(X)$ its quotient by the subgroup of classes $F$ such that $\chi(E,F)=0$ for all $E\in D(X)$, where $\chi(E,F)$ is the Euler form
$$ \chi(E,F) = \sum_i \dim\Hom^i_{D(X)}(E,F). $$

\begin{defi}
A \textit{full numerical stability condition} on $X$ is a pair $\sigma=(Z,\cA)$ whereby
\begin{itemize}
	\item $\cA$ is the \textit{heart} of a bounded $t$-structure on $D(X)$.
	\item $\cZ:K_{num}(X)\fun\bC$ a group homomorphism called the \textit{central charge}.
\end{itemize}
satisfying properties 1,2 and 3 below.
\begin{description}
	\item[1 (Positivity)] For all $0\neq E\in\cA$, $\cZ(E)\in\{re^{i\pi\varphi}\mid r>0,0<\varphi\leq1\}$.
\end{description}

To discuss stability for a given stability condition, we define for each $E\in D(X)$
$$ \beta(E) = -\frac{\Real \cZ(E)}{\Im \cZ(E)} \in (-\infty, \infty] $$
For example, if $\cZ(E)=-1$, then $\beta(E) = \infty$, and if $\cZ(E)=\sqrt{-1}$, then $\beta(E) = 0$. 

We say that $E\in\cA$ is $\sigma$-stable (resp.\ $\sigma$-semistable) if for all nontrivial $F\su E$ in $\cA$, we have $\beta(E)>\beta(F)$ (resp.\ $\beta(E)\geq\beta(F)$).

\begin{description}
	\item[2 (Harder-Narasimhan Filtrations)] For all $E\in\cA$, there exist objects $E_1,\ldots,E_{n-1}\in\cA$ such that
		\begin{itemize}
			\item $0=E_0\su E_1\su\cdots\su E_{n-1}\su E_n=E$ in $\cA$.
			\item $E_{i+1}/E_i$ is $\sigma$-semistable for each $i$.
			\item $\beta(E_{1}/E_0)>\beta(E_{2}/E_1)>\cdots>\beta(E_{n}/E_{n-1})$.
		\end{itemize}
		
	\item[3 (Support Property)] Choose a norm $\|.\|$ on $K_{num}(X)\otimes\bR$. There exists a $C>0$ such that for all $\sigma$-semistable $E\in D(X)$, we have $C\|E\|\leq|\cZ(E)|$.
\end{description}
\end{defi}

\begin{remark}
The support property guarantees a nicely behaved wall and chamber structure for classes of objects - namely, the walls are locally finite, real codimension 1 submanifolds of the stability manifold, and deleting the walls gives chambers where Bridgeland stability is constant (see \cite[Proposition 3.3]{stabcondsonlocalp2}). The support property is equivalent to Bridgeland's concept of \textit{full} (see \cite[Proposition B.4]{stabcondsonlocalp2}).
\end{remark}

Throughout the rest of this paper, when we say \textbf{stability condition}, we mean a full numerical stability condition.


\subsection{Bridgeland stability conditions on a surface}\label{stabcond}

Let $S$ be a smooth projective surface.
The stability conditions that we are going to consider were defined in \cite{ABL}.
They form a subset $\stabdivS$ of stability conditions that depend on a choice of an ample and a general divisor (the ``div'' stands for ``divisor'') and are well-suited to computations.
Let us recall their definition.

Let $S$ be a smooth projective surface. Given two $\bR$-divisors $D,H$ with $H$ ample, we define a stability condition $\sigma_{D,H}=(\cZ_{D,H},\cA_{D,H})$ on $S$ as follows:

Consider the $H$-Mumford slope
$$ \mu_H(E) = \ds\frac{c_1(E).H}{\rk(E)H^2}. $$

Let $\cA_{D,H}$ be the tilt of the standard $t$-structure on $D(S)$ at $\mu_H(D) =\ds\frac{D.H}{H^2}$ defined by $\cA_{D,H} = \{E\in D(S) \mid H^i(E)=0 \text{ for } i\neq-1,0,\ H^{-1}(E)\in\cF_{D,H},\ H^0(E)\in\cT_{D,H}\}$, where
\begin{itemize}
\item
$\cT_{D,H}\su\Coh(S)$ is the full subcategory closed under extensions generated by torsion sheaves and $\mu_H$-stable sheaves $E$ with $\mu_H(E)>\ds\frac{D.H}{H^2}$.
\item
$\cF_{D,H}\su\Coh(S)$ is the full subcategory closed under extensions generated by $\mu_H$-stable sheaves $F$ with $\mu_H(F)\leq\ds\frac{D.H}{H^2}$.
\end{itemize}
\vspace{.1in}

Now, define $\cZ_{D,H}$ as
$ \cZ_{D,H}(E) = - \ds\int{e^{-(D+iH)}\ch(E)}. $
This is equal to
$$ 
\left( -\ch_2(E) + c_1(E).D - \frac{\rk(E)}{2}(D^2 - H^2) \right)+ i \left( c_1(E).H - \rk(E) D.H \right). $$

According to \cite[Corollary 2.1]{ABL} and \cite[Sections 3.6 \& 3.7]{stabcondsextrcontrs}, $\sigma_{D,H}$ is a stability condition on $S$. Let $\stabdivS$ be the set of all such stability conditions.

\begin{remark}
These are \textit{geometric} stability conditions because for all $p\in S$, the skyscraper sheaf $\bC_p \in \cA_{D,H}$ is $\sigma_{D,H}$-stable with $\cZ_{D,H}(\bC_p)=-1$ (the proof is the same as \cite[Proposition 6.2.a]{ABCH}).
\end{remark}

{\it Note}.
When the divisors $D$ and $H$ have been fixed, we will often drop the $D,H$ subscript from $\sigma$, $\cZ$, $\cA$, $\cT$, and $\cF$.

\subsection{Slices of $\stabdivS$}\label{slices}

Let $H$ be an ample $\bR$-divisor such that $H^2=1$.
If $S$ has Picard rank 1, then the stability conditions in $\stabdivS$ are all of the form $\sigma_{sH,tH}$.
It was already proved in \cite{ABL} that line bundles are always Bridgeland stable for these stability conditions.
Assume herein that $S$ has Picard rank greater than 1.

One of the features that makes $\stabdivS$ well suited to computations is its decomposition into well-behaved 3-spaces, each given by a choice of ample divisor and another divisor orthogonal to it.
In these 3-spaces, walls of interest will be quadric surfaces, and most of our work will begin by first choosing a particular 3-space to live in.
Most of these concepts were introduced in \cite{MaciociaWalls}.

\begin{defi}
Choose an $\bR$-divisor $G$ with $G.H=0$ and $G^2=-1$ (note that $G^2\leq0$ by the Hodge Index Theorem with $G^2=0$ iff $G=0$).
Then, define $\cS_{G,H} := \{\sigma_{sH+uG,tH}\mid s,u,t\in\bR,t>0\} \su \stabdivS$.
\end{defi}

Herein, we assume that any divisors $G,H$ are as above. We identify $\cS_{G,H}$ with $\{ (s,u,t)\ |\ t>0 \}$ by $(s,u,t) \leftrightarrow \sigma_{sH+uG,tH}$.  

Each of the stability conditions $\sigma_{D,H}$ defined in \ref{stabcond} can be observed as an element of a particular 3-dimensional slice. Indeed, we can simply scale $H$ to ensure that $H^2=1$ and then choose $G$ such that $D=sH+uG$ so that $\sigma_{D,H} \in \cS_{G,H}$. Thus these slices cover all of $\stabdivS$.

{\it Note}. Although the spaces $\cS_{G,H}$ do not contain the plane $t=0$, we will, for convenience, expand them to include the plane $t=0$.
This is because when studying a potential wall as an equation in $s,u,t$, it will become useful to look at the intersection of the wall with the plane $t=0$ (see Section \ref{walls} below for more details).
In addition, we will 
write $(s,u)$ to mean $(s,u,0)$.

Let $E\in\DS$, and set $\ch(E)=(\ch_0(E),\ch_1(E),\ch_2(E)) = (r,c_1(E),c)$. We may write $c_1(E)=d_hH+d_gG+\alpha$, where $\alpha.H=\alpha.G=0$ and $d_h,d_g\in\bR$.
Specifically, we have $d_h=c_1(E).H$ and $d_g=-c_1(E).G$.

\begin{remark}\label{ssVerticalPlane}
Because $\mu_{H}(sH+uG)=s$, the vertical plane $s=\mu_H(E)$ plays a very important role in $\cS_{G,H}$.
Indeed, if $E$ is a $\mu_H$-semistable sheaf, then $s<\mu_H(E)$ iff $E\in\cA_{sH+uG,H}$, and $s\geq\mu_H(E)$ iff $E[1]\in\cA_{sH+uG,H}$.
\end{remark}

\begin{defi}
For each fixed value of $u$, we denote by $\Pi_u$ the vertical plane of stability conditions $(s,u,t)$ of fixed $u$-value.
It is parametrized by $(s,t)$, and it will play a special role in our work (see Section \ref{walls} below for more details).
\end{defi}

The central charge of a stability condition $\sigma=\sigma_{sH+uG,tH}$ for an object $E$ with $\ch(E)=(r,d_hH+d_gG+\alpha,c)$ is equal to
$$ \cZ(E) = \left( -c + s d_h - u d_g - \frac{r}2(s^2 - u^2 - t^2) \right)+ i \left( t d_h - r s t \right). $$


\section{Reduction to the case of $\cO_S$}
\label{reduction}

The action of tensoring stability conditions by line bundles will allow us to restrict our attention from the stability of all line bundles to that of $\cO_S,\cO_S[1]$.

We prove the following statement:
\begin{lem}\label{reductionlemma}
Let $L=\cO_S(D_1)$ be a line bundle, and let $\sigma_{D,H}$ be a Bridgeland stability condition as above.
Then, $E$ destabilizes $L$ $($resp.\ $L[1])$ at $\sigma_{D,H}$ if and only if $\cO_S(-D_1)\otimes E$ destabilizes $\cO_S$ $($resp.\ $\cO_S[1])$ at $\sigma_{D-D_1,H}$.
\end{lem}

{\bf Notation}: We shall denote $\sigma_{D-D_1,H}$ by $\cO_S(-D_1)\otimes\sigma_{D,H}$.

\vspace{.1cm}

\pf
We show the proof for $L$ (the case of $L[1]$ is similar).
To simplify the notation within the proof, let $\sigma=\sigma_{D,H}$ and $\sigma'=\sigma_{D-D_1,H}$.

If $E$ destabilizes $L$ at $\sigma$, then we must have that $E\su L$ in $\cA_{D,H}$ and that $\beta_{\sigma}(E)>\beta_{\sigma}(L)$.
Let us prove that, in this case, $\cO_S(-D_1)\otimes E\su\cO_S$ in $\cA_{D-D_1,H}$ and that $\beta_{\sigma'}(\cO_S(-D_1)\otimes E)>\beta_{\sigma'}(\cO_S)$.
It suffices to prove the following statements about an object $E'\in\cA_{D,H}$:
\begin{itemize}
\item
$E'\in\cA_{D,H}$ iff $\cO_S(-D_1)\otimes E'\in\cA_{D-D_1,H}$.
\item
$\cZ_{D,H}(E')=\cZ_{D-D_1,H}(\cO_S(-D_1)\otimes E')$.
\end{itemize}

Recall that $\cA_{D,H}$ is the full subcategory of $D(S)$ closed under extensions generated by torsion sheaves, $\mu_H$-stable sheaves with $\mu_H>\mu_H(D)$, and shifts of $\mu_H$-stable sheaves with $\mu_H\leq\mu_H(D)$.
Because, given a $\mu_H$-stable sheaf $E''$, we have that $\mu_H(E'')>\mu_H(D)$ iff $\mu_H(\cO_S(-D_1)\otimes E'')>\mu_H(D-D_1)$, we can conclude that an object $E'\in\cA_{D,H}$ iff $\cO_S(-D_1)\otimes E'\in\cA_{D-D_1,H}$.

As for the second statement,
\begin{align*}
\cZ_{D,H}(E')
&= - \int{e^{-(D+iH)}\ch(E')} \\
&= - \int{e^{-(D-D_1+iH)}e^{-D_1}\ch(E')} \\
&= - \int{e^{-(D-D_1+iH)} \left( 1[S] - D_1 + \frac{(-D_1)^2}{2}[pt] \right).\ch(E')} \\
&= - \int{e^{-(D-D_1+iH)} \ch(\cO_S(-D_1)).\ch(E')} \\
&= - \int{e^{-(D-D_1+iH)} \ch(\cO_S(-D_1) \otimes E')} \\
&= \cZ_{D-D_1,H}(\cO_S(-D_1) \otimes E').
\end{align*}

The converse (i.e., if $\cO_S(-D_1)\otimes E$ destabilizes $\cO_S$ at $\sigma'$, then $E$ destabilizes $L$ at $\sigma$) is also true because we can start with $E\otimes\cO_S(-D_1)\su\cO_S$ in $\cA_{D-D_1,H}$ and tensor through by $\cO_S(D_1)$.
\qed


\section{Preliminaries on the stability of $\cO_S$}
\label{prelims}

Walls are subsets of the stability manifold where the stability of objects can change. Our main interest lies in describing the chambers of stability for $\cO_S$, which are bounded by walls corresponding to certain destabilizing objects.
Let us start with a few definitions and a description of the possible walls.


\subsection{Subobjects of $\cO_S$ and their walls}\label{walls}

First, here is our generic definition of a wall.

\begin{defi}
Given two objects $E,B \in \DS$, with $B$ Bridgeland-stable for at least one stability condition, we define the wall $\cW(E,B)$ as $\{ \sigma\in\stabdivS\ |\ (\Real \cZ(E))(\Im \cZ(B))-(\Real \cZ(B))(\Im \cZ(E))=0 \}$. 
If at some $\sigma\in\cW(E,B)$ we have $E\su B$ in $\cA$, we say that $\cW(E,B)$ is a \textbf{weakly destabilizing wall} for $B$.
If at some $\sigma\in\cW(E,B)$ we have $E\su B$ in $\cA$, and $B$ is Bridgeland $\sigma$-semistable, we say that $\cW(E,B)$ is an \textbf{actually destabilizing wall} for $B$. 
\end{defi}

Note that if $\Im \cZ(E)\neq0\neq \Im \cZ(B)$ then the defining condition is simply $\beta(E)=\beta(B)$.
We are interested in the walls for $\cO_S$ and $\cO_S[1]$, and we start by studying the walls for $\cO_S$.
Let us assume that $\cO_S\in\cA$.
At each fixed value of $u$, Maciocia showed in \cite[Section 2]{MaciociaWalls} that all walls for $\cO_S$ in $\Pi_u$ are nested semicircles centered on the $s$-axis. Therefore, given two objects $E_1$ and $E_2$ and a fixed value of $u$, we have that $\cW(E_1,\cO_S)\cap\Pi_u$ and $\cW(E_2,\cO_S)\cap\Pi_u$ are both semicircles, with one of them inside the other, unless they are equal.

\begin{figure}[h]
\begin{center}
\includegraphics[scale=.4]{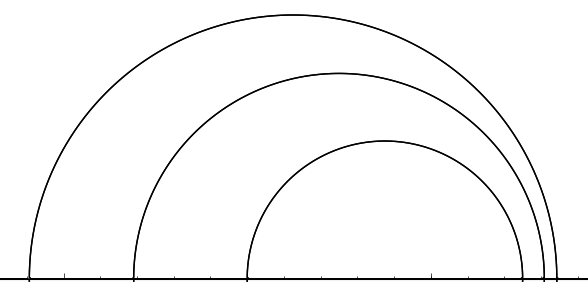}
\label{nested}
\caption{Nested walls in $\Pi_u$}
\end{center}
\end{figure}

\begin{defi}
We say that the wall $\cW(E_1,\cO_S)$ is {\bf inside} the wall $\cW(E_2,\cO_S)$ at $u$ if the semicircle $\cW(E_1,\cO_S)\cap\Pi_u$ is inside the semicircle $\cW(E_2,\cO_S)\cap\Pi_u$ or equal to it.
We will use the notation
$$ \cW(E_1,\cO_S)\cap\Pi_u \preceq \cW(E_2,\cO_S)\cap\Pi_u. $$
When this happens, we will also refer to the second wall as the {\bf higher} wall at $u$.
\end{defi}

\begin{lem}
\label{subsofO}
Let $\sigma\in\stabdivS$, and let $0\fun E\fun\cO_S\fun Q\fun0$ be a short exact sequence in $\cA$.
Then, $E$ is a torsion-free sheaf, and $H^0(Q)$ is a quotient of $\cO_S$ of rank $0$.
In particular, the kernel of the map $\cO_S\fun H^0(Q)$ is an ideal sheaf $\cI_Z(-C)$ for some effective curve $C$ and some zero-dimensional scheme $Z$ $($with $C$ or $Z$ possibly $0)$.
\end{lem}

\pf
The long exact sequence in cohomology associated to the short exact sequence shows that $E$ must be a sheaf, and $Q$ may have cohomologies of degrees $-1$ and $0$:
$$ 0 \lfun H^{-1}(Q) \lfun E \lfun \cO_S \lfun H^0(Q) \lfun 0. $$
If $H^0(Q)$ is of rank 1, then it would have to be equal to $\cO_S$, and we would have that $H^{-1}(Q)\simeq E$.
However, $H^{-1}(Q)\in\cF$, $E\in\cT$, and $\cF\cap\cT=\{0\}$.
Therefore, in this case, both $H^{-1}(Q)$ and $E$ would have to be $0$.
This means that $H^0(Q)$ must be a quotient of $\cO_S$ of rank $0$.
The kernel of $\cO_S\fun H^0(Q)$ is a torsion-free sheaf of the form $\cI_Z(-C)$ for some effective curve $C$ and some zero-dimensional scheme $Z$, with $C$ or $Z$ possibly $0$.
Because $E$ is an extension of the torsion-free sheaves $\cI_Z(-C)$ and $H^{-1}(Q)$, it is also a torsion-free sheaf.
\qed

\vspace{.1in}

Here, we study which forms the walls $\cW(E,\cO_S)$ can take in $\cS_{G,H}$.
The intersection of a wall with the $t=0$ plane is a conic through the origin, and we classify the wall based on invariants associated to $E$.

If $\ch(E)=(r,d_hH+d_gG+\alpha_E,c)$, then we observed above that
$$ \cZ(E)=\left( -c + s d_h - u d_g - \frac{r}2(s^2 - u^2 - t^2) \right)+ i \left( t d_h - r s t \right), $$
and the equation of the wall $\cW(E,\cO_S)$ is
$$ \frac{t}2 (-d_h (s^2 + t^2 + u^2) + 2 d_g s u + 2 c s) = 0. $$
Because $t\neq0$, this is equivalent to
$$ -d_h (s^2 + t^2 + u^2) + 2 d_g s u + 2 c s = 0. $$
We now think of the $\cS_{G,H}$ spaces as being extended to the $t=0$ plane, and we study these quadrics by studying their intersection with the $t=0$ plane:
$$ - d_h (s^2 + u^2) + 2 d_g s u + 2 c s = 0. $$

Note that because the walls are semicircles in $\Pi_u$ for any fixed $u$ by \cite{MaciociaWalls}, knowing where the wall is at $t=0$ would tell us where the wall is at any $t>0$.
This is why it is sufficient to study the quadrics in the $t=0$ plane.

We will abuse notation and continue to refer to the $t=0$ equation above as the wall $\cW(E,\cO_S)$.
Its discriminant is equal to
$$ \Delta = 4 (d_g^2 - d_h^2). $$

If $\Delta=0$, then the wall is a parabola.
Because $\cO_S\in\cA$, we have that $s<0$ by Remark \ref{ssVerticalPlane}.
Therefore, the wall can only be a weakly destabilizing wall if $c>0$, in which case the equation of the parabola is $-d_h(s\pm u)^2+2cs=0$.

Assume now that $\Delta\neq0$.
We can show the following:

\begin{lem}\label{PW}
Assume that $\Delta\neq0$ and $c\neq0$.
Then,
\begin{itemize}
\item
The following two points are on the wall, and the tangent line to these points is vertical (i.e., of the form $s=\,$constant):
$$ (0,0) \textrm{ and } P_W := -\frac{2c}{d_g^2 - d_h^2} (d_h, d_g). $$
\item
The tangent line to the wall is horizontal (i.e., of the form $u=\,$constant) at the points where the conic intersects $u=s$ and $u=-s$, which are
$$ \left( \frac{c}{d_h-d_g}, \frac{c}{d_h-d_g} \right) \textrm{ and } 
\left( \frac{c}{d_h+d_g}, -\frac{c}{d_h+d_g} \right). $$
\end{itemize}
\end{lem}

\pf
We only show the first result. The proof of the second one is very similar.

To check that the points are on the wall, it suffices to check that $f(0,0)=f(P_W)=0$.
Clearly, $f(0,0)=0$.
If we plug in $P_W$, we obtain
$$ f(P_W) = - d_h \cdot \frac{4c^2}{(d_g^2-d_h^2)^2} (d_h^2 + d_g^2) + 2 d_g \cdot \frac{4c^2}{(d_g^2-d_h^2)^2} d_h d_g - \frac{4c^2}{d_g^2-d_h^2} d_h. $$
Because $\Delta\neq0$ and $c\neq0$, we can multiply through by $(d_g^2-d_h^2)/(4c^2)\neq0$ and obtain
$$ - d_h (d_h^2 + d_g^2) + 2 d_g d_h d_g - (d_g^2-d_h^2) d_h = -d_h^3 - d_h d_g^2 + 2 d_h d_g^2 - d_h d_g^2 + d_h^3, $$ 
which is equal to $0$.

To check that the tangent line at those points is vertical, it suffices to prove that the derivative of $f(s,u)$ with respect to $u$ is $0$ there.

We have that
$$ f_u(s,u) = - 2 d_h u + 2 d_g s. $$
Then, $f_u(0,0)=0$, and
$$ f_u(P_W) = \frac{4c}{d_g^2-d_h^2}d_hd_g - \frac{4c}{d_g^2-d_h^2}d_gd_h = 0. $$
\qed

If $\Delta<0$, then the wall is an ellipse.
It can only be a weakly destabilizing wall if the ellipse is to the left of $s=0$, and this only occurs if $c>0$ (so that $P_W$ is to the left of $s=0$).

If $\Delta>0$, then there are three possibilities:
\begin{itemize}
\item
If $c=0$, then the wall is a cone centered at $(0,0)$.
\item
If $c>0$, then the wall is a hyperbola with center $P$ to the right of $s=0$.
\item
If $c<0$, then the wall is a hyperbola with center $P$ to the left of $s=0$.
In this case, the asymptotes have slope
$$ \frac{d_g\pm\sqrt{d_g^2-d_h^2}}{d_h}. $$
\end{itemize}

Below is a summary of all possible weakly destabilizing walls (pictures drawn for $d_h<0$ and $d_g>0$):
\begin{description}
\item[Parabola.]
When $d_g^2-d_h^2=0$ and $c>0$.
\item[Ellipse.]
When $d_g^2-d_h^2<0$ and $c>0$.
\item[Cone.]
When $d_g^2-d_h^2>0$ and $c=0$.
\item[Right Hyperbola.]
When $d_g^2-d_h^2>0$ and $c>0$.
\item[Left Hyperbola.]
When $d_g^2-d_h^2>0$ and $c<0$.
\end{description}

\begin{figure}[h]
\centering
\begin{minipage}{.32\textwidth}
  \centering
  \includegraphics[scale=.15]{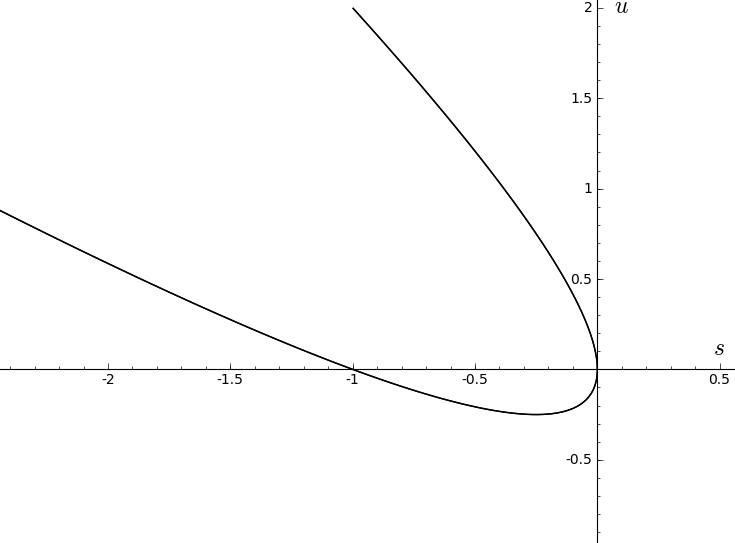}
  \label{Parwall}
  \caption{Parabola} 
\end{minipage}%
\begin{minipage}{.32\textwidth}
  \centering
  \includegraphics[scale=.15]{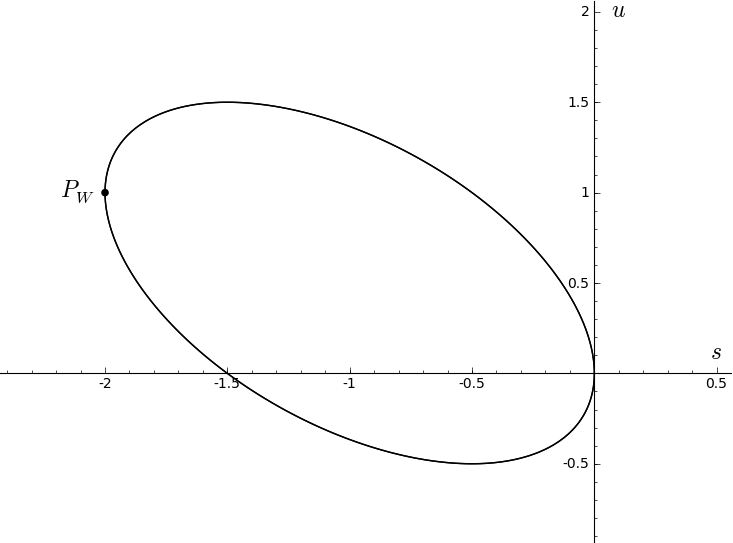}
  \caption{Ellipse}
  \label{Ellwall}
\end{minipage}
\begin{minipage}{.32\textwidth}
  \centering
  \includegraphics[scale=.15]{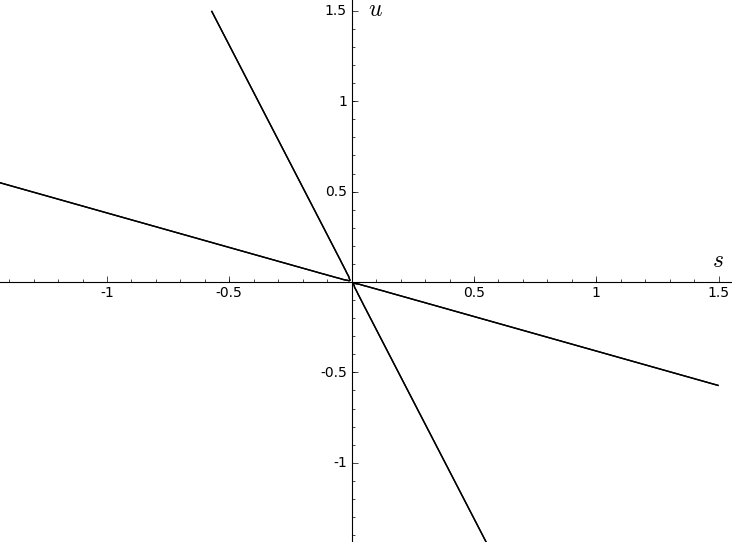}
  \caption{Cone}
  \label{Conwall}
\end{minipage}
\end{figure}

\begin{figure}[h]
\centering
\begin{minipage}{.5\textwidth}
  \centering
  \includegraphics[scale=.21]{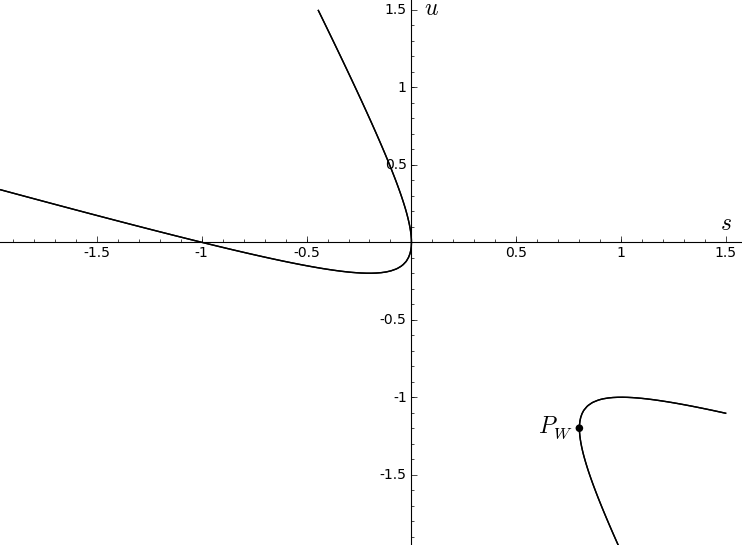}
  \label{HypRWall}
  \caption{Right hyperbola} 
\end{minipage}%
\begin{minipage}{.5\textwidth}
  \centering
  \includegraphics[scale=.21]{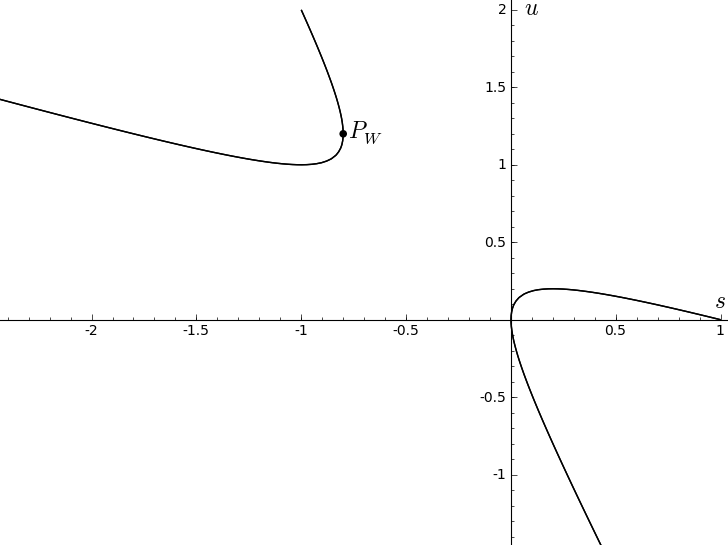}
  \caption{Left hyperbola}
  \label{HypLWall}
\end{minipage}
\end{figure}


We end this section by noting an important geometric property of these walls that will be useful in various proofs later in the paper.

\begin{lem}\label{geometry}
Given two subobjects $E_1$ and $E_2$ of $\cO_S$ in $\cA$, there exists at most one non-zero value $u_0$ of $u$ such that
$$ \cW(E_1,\cO_S)\cap\Pi_{u_0} = \cW(E_2,\cO_S)\cap\Pi_{u_0}, $$
unless the two walls coincide everywhere.
\end{lem}

\pf
Suppose that $ \cW(E_1,\cO_S) $ and $ \cW(E_2,\cO_S) $ do not coincide everywhere.
Looking at the intersection of the walls with the $t=0$ plane, we see that they are two conics that intersect at $(0,0)$ with multiplicity two.
Therefore, they can only intersect at at most two other points in the $t=0$ plane.
Because the walls $\cW(E,\cO_S)\cap\Pi_u$ are nested semicircles centered on the $s$-axis, they intersect the $t=0$ plane.
This means that, if $\cW(E_1,\cO_S)\cap\Pi_{u_0} = \cW(E_2,\cO_S)\cap\Pi_{u_0}$ for some $u_0$, the walls coincide at the two points where the semicircle intersects the $t=0$ plane and cannot intersect anywhere else except at $(0,0)$, thus proving that the walls cannot coincide at any other non-zero value of $u$.
\qed



\subsection{Bridgeland stability of $\cO_S$ for $t>>0$}

We now note a statement that is already known and will be needed below.

\begin{pro}\label{Olarget}
Let $H$ and $G$ be as above, and fix a divisor $sH+uG$ with $s<0$.
Then, $\cO_S$ is Bridgeland stable for the stability condition $\sigma_{sH+uG,tH}$ for $t>>0$.
\end{pro}

\pf
This follows, for example, from the result in \cite[Theorem 3.11]{MaciociaWalls} that walls in each plane $\Pi_u$ are disjoint semicircles that are bounded above.
\qed

\begin{remark}
Even when fixing $G$ and $H$, a $t_0$ such that $\cO_S$ is Bridgeland stable for all $(s,u,t)$ with $t>t_0$ does not necessarily exist.
We will see below that whenever $\cO_S$ is destabilized by a line bundle $\cO_S(-C)$ in a space $\cS_{G,H}$, the destabilizing wall will be a hyperboloid, and there will be stability conditions for any $t$ for which $\cO_S$ is not Bridgeland stable.
However, for every fixed value of $u$, Maciocia proved in \cite{MaciociaWalls}  that there exists a value $t_0(u)$ such that $\cO_S$ is Bridgeland stable for $\sigma_{sH+uG,tH}$ for all $s<0$ and $t>t_0(u)$.
\end{remark}


\subsection{Bertram's Lemma}\label{bertram}

The following lemma is a key tool in the proof of the main theorem. The lemma is essentially Lemma 6.3 in \cite{ABCH} adapted to our situation.
It allows us to, in some situations, find walls that are higher than a given wall for $\cO_S$ by removing a Mumford semistable factor from the subobject or quotient. In either case, the rank strictly decreases, thus providing the foundation for our induction proof characterizing the stability of $\cO_S$.

Let $E$ be a subobject of $\cO_S$ of rank $\geq2$ in $\cA$, and let $Q$ be the quotient.

Let $0=E_0\su E_1\su E_2\su\cdots\su E_{n-1}\su E_n=E$ be the Harder-Narasimhan filtration of $E$, and let $K_i=E_i/E_{i-1}$ (so that $K_1=E_1$ and $K_n=E/E_{n-1}$).
We have that $\mu_H(K_1)>\mu_H(K_2)>\cdots>\mu_H(K_n)$.
In addition, let $\Kl=K_n$.

Similarly, let $0=F_0\su F_1\su F_2\su\cdots\su F_{m-1}\su F_m=H^{-1}(Q)$ be the Harder-Narasimhan filtration of $H^{-1}(Q)$, and let $J_i=F_i/F_{i-1}$ with $\Jl=J_1$.

\begin{remark}\label{In_cA}
Note that $E\su\cO_S\in\cA$ iff $\mu_H(\Jl)\leq s<\mu_H(\Kl)$.
Moreover, we must have that $\mu_H(\Jl)<\mu_H(\Kl)<0$.
\end{remark}

\pf
For $E$ to be a subobject of $\cO_S$ in $\cA$, we need in particular that both $E$ and $H^{-1}(Q)$ be objects of $\cA$.
For $E$ to be an object of $\cA$, we must have that $s<\mu_H(\Kl)$, and for $H^{-1}(Q)$ to be an object of $\cA$, we must have that $s\geq\mu_H(\Jl)$.
The last condition that $\mu_H(\Kl)<0$ originates from the fact that, as we observed above, $E$ is an extension of $\cI_Z(-C)$ by $H^{-1}(Q)$ for some effective curve $C$ and zero-dimensional scheme $Z$ (with $C$ or $Z$ possibly zero).
Because $\cO_S\in\cA$ implies that $s<0$ by Remark \ref{ssVerticalPlane}, then $\mu_H(\Jl)<s<0$ implies that $\mu_H(H^{-1}(Q))\leq\mu_H(\Jl)<0$.
Because $E$ is an extension of $\cI_Z(-C)$ by $H^{-1}(Q)$, we must also have that $\mu_H(E)<0$, and therefore, $\mu_H(\Kl)\leq\mu_H(E)<0$.
\qed

The idea for Bertram's lemma is the following: For every value of $u$, the wall $\cW(E,\cO_S)\cap\Pi_u$ is a semicircle in the upper-half plane $\Pi_u$ parametrized by $(s,t)$.
The lemma shows that if the semicircle intersects the lines $s=\mu_H(\Jl)$ or $s=\mu_H(\Kl)$, then  there exists another subobject $E'\su\cO_S$ such that $\cW(E,\cO_S)$ is inside $\cW(E',\cO_S)$ at $u$.
If the semicircle $\cW(E,\cO_S)\cap\Pi_u$ intersects the line $s=\mu_H(\Kl)$, then $E'=E_{n-1}$, and if it intersects the line $s=\mu_H(\Jl)$, then $E'=E/\Jl$ (note that, as sheaves, $\Jl=J_1=F_1\su H^{-1}(Q)\su E$).


\begin{lem}[Bertram's Lemma]\label{BertramLemma}
Fix $H$ and $G$ as above, and let $E\su\cO_S$ in $\cA_0$ for some $\sigma_0=\sigma_{s_0H+u_0G,t_0H}=(\cZ,\cA_0)$ such that $\sigma_0\in\cW(E,\cO_S)$.
\begin{itemize}
\item[$(1)$]
If $\cW(E,\cO_S)\cap\Pi_{u_0}$ intersects the line $s=\mu_H(\Kl)$ for $t>0$, then $\beta(E_{n-1})>\beta(E)$ at $\sigma_0$, with $E_{n-1}\su\cO_S$ in $\cA_0$ $($in particular, $E$ is not $\mu_H$-semistable$)$.
\item[$(2)$]
If $\cW(E,\cO_S)\cap\Pi_{u_0}$ intersects the line $s=\mu_H(\Jl)$ for $t>0$, then $\beta(E/\Jl)>\beta(E)$ at $\sigma_0$, with $E/\Jl\su\cO_S$ in $\cA_0$.
\end{itemize}
\end{lem}

\pf
(1) Maciocia proved in \cite{MaciociaWalls} that if $E$ is $\mu_H$-semistable, then $\cW(E,\cO_S)\cap\Pi_u$ does not intersect the line $s=\mu_H(E)$ for any $u$.
Therefore, if $\cW(E,\cO_S)\cap\Pi_{u_0}$ intersects $s=\mu_H(\Kl)$, $E$ cannot be $\mu_H$-semistable.
Because $E\in\cA_0$, we have that $s_0<\mu_H(\Kl)$.
Because $\Kl\in\cA$ iff $s<\mu_H(\Kl)$, it follows that for all values of $s$ between $s_0\leq s<\mu_H(\Kl)$, there is a short exact sequence $0\fun E_{n-1}\fun E\fun\Kl\fun0$ in $\cA$.
At $s=\mu_H(\Kl)$, we have $\Im \cZ(\Kl)=0$, and $\beta(\Kl)=-\infty$.
Therefore, approaching $s=\mu_H(\Kl)$ from $\sigma_0$ along $\cW(E,\cO_S)\cap\Pi_{u_0}$, we must have that $\beta(\Kl)\fun-\infty$, and $\beta(E_{n-1})>\beta(E)>\beta(\Kl)$.
Because the walls are nested in $\Pi_{u_0}$, it follows that $\beta(E_{n-1})>\beta(E)=\beta(\cO_S)$ at $\sigma_0$ as well.

\begin{figure}[h]
\centering
\begin{minipage}{.5\textwidth}
  \centering
  \includegraphics[scale=.21]{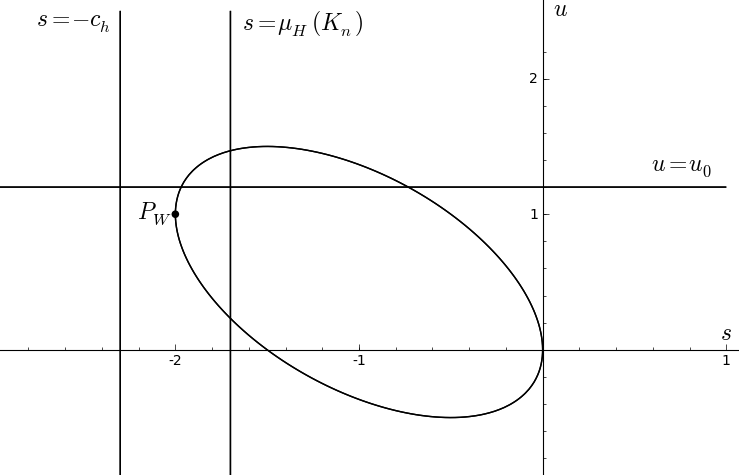}
  \caption{$t=0$ plane} 
\end{minipage}%
\begin{minipage}{.5\textwidth}
  \centering
  \includegraphics[scale=.21]{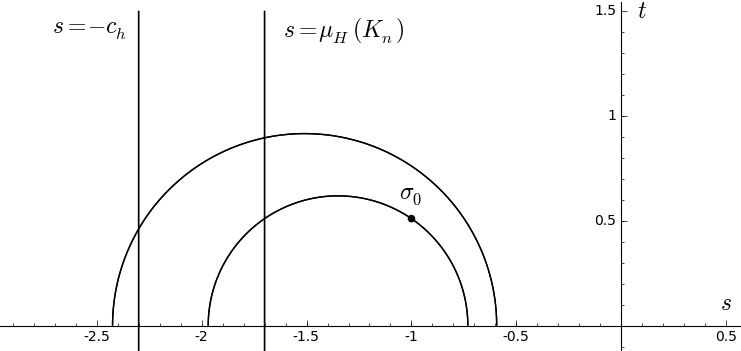}
  \caption{$\Pi_u$ plane at $u=u_0$}
\end{minipage}
\end{figure}

(2) Because $E\su\cO_S\in\cA_0$, $s_0$ must satisfy the following inequalities: $\mu_H(\Jl)\leq s_0<\mu_H(\Kl)$.
Therefore, for all values of $s$ between $\mu_H(\Jl)\leq s\leq s_0$, we have that $E,\Jl[1]\in\cA$, and there exists a short exact sequence $0\fun E\fun E/\Jl\fun\Jl[1]\fun0$ in $\cA$.

At $s=\mu_H(\Jl)$, we know that $\Im \cZ(\Jl[1])=0$, and $\beta(\Jl[1])=\infty$.
If we approach $s=\mu_H(\Jl)$ from $\sigma_0$ along $\cW(E,\cO_S)\cap\Pi_{u_0}$, we must have that $\beta(\Jl[1])\fun\infty$, and $\beta(E)<\beta(E/\Jl)<\beta(\Jl[1])$.
Because the walls are nested in $\Pi_{u_0}$, it follows that $\beta(E/\Jl)>\beta(E)=\beta(\cO_S)$ at $\sigma_0$ as well.
\qed

\begin{remark}
Because we are considering surfaces with Picard rank greater than 1, we are forced to strengthen the hypotheses from those of \cite{ABCH}, but the proof is identical. We have restricted our attention to $\cO_S$, but the proof holds for any sheaf satisfying the Bogomolov inequality (even a torsion sheaf).
\end{remark}

Let us note an important fact that will be needed later in the paper:

\begin{remark}
\label{emodjremark}
Let $E,Q,\Jl$ be as above, with $E\su\cO_S$ in $\cA$, $Q$ the quotient of $\cO_S$ by $E$ in $\cA$, and $\Jl$ the first non-zero term in the Harder-Narasimhan filtration of $H^{-1}(Q)$.
Now, consider $E/\Jl\su\cO_S$ in $\cA$, and let $Q'$ be the quotient of $\cO_S$ by $E/\Jl$ in $\cA$.
Then, we have that $H^{-1}(Q')=H^{-1}(Q)/\Jl$ and $H^0(Q')=H^0(Q)$.
\end{remark}


\section{Bridgeland stability of $\cO_S$}
\label{stabofO}

In this section, we prove the first part of Conjecture \ref{mainconj} for surfaces with any Picard rank and no curves of negative self-intersection as well as surfaces with Picard rank 2 and one irreducible curve of negative self-intersection.
Proposition \ref{conjevidencePicRk2} serves as evidence for the conjecture on surfaces with Picard rank 2 and two irreducible curves of negative self-intersection.

\subsection{Subobjects of $\cO_S$ of rank $1$}

Understanding the rank-1 weakly destabilizing subobjects of $\cO_S$ is crucial to our main results, all of which use induction. If a subobject $E\su\cO_S$ is of rank $1$, then Lemma \ref{subsofO} shows that $E$ must be equal to $\cI_Z(-C)$ for some effective curve $C$ and some zero-dimensional scheme $Z$, with $C$ or $Z$ possibly $0$.

If $C=0$, then $E=\cI_Z$ for some zero-dimensional scheme $Z$.
In this case, $\cI_Z$ does not destabilize $\cO_S$ because $\beta(\cI_Z)=(-2l(Z)+s^2-u^2-t^2)/(-2st),$ where $l(Z)$ is the length of $Z$, and $\beta(\cO_S)=(s^2-u^2-t^2)/(-2st)$.
Therefore, when $\cO_S\in\cA$, we have that $s<0$, and $\beta(\cI_Z)<\beta(\cO_S)$.
This means that $\cI_Z$ does not destabilize $\cO_S$ whenever $\cO_S\in\cA$.

Assume now that $C\neq0$, and consider the subobject $\cI_Z(-C)\su\cO_S$, with $Z$ a zero-dimensional scheme $Z$ (possibly equal to $0$).
We show that for $\cI_Z(-C)$ to weakly destabilize $\cO_S$, we must have $C^2<0$.


\begin{pro}\label{rk1}
If $C^2\geq0$, then $\cI_Z(-C)$ does not weakly destabilize $\cO_S$ for any $\sigma$, i.e., there does not exist any $\sigma$ such that $\cI_Z(-C)\su\cO_S$ in $\cA$ and $\beta(\cI_Z(-C))=\beta(\cO_S)$.
\end{pro}

\pf
Note that $\ch(\cI_Z(-C))=(1,-C,C^2/2-l(Z))$, where $l(Z)$ is the length of $Z$.
If $C=c_hH+c_gG+\alpha_C$, with $\alpha_C.H=\alpha_C.G=0$, then $C^2=c_h^2-c_g^2+\alpha_C^2$, and $\alpha_C^2\leq0$ by the Hodge Index Theorem.
Therefore, if $C^2\geq0$, then $c_h^2-c_g^2\geq0$.
The equation for the wall $\cW(\cI_Z(-C),\cO_S)$ simplifies to
$$ c_h (s^2 + t^2 + u^2) - 2 c_g s u + (c_h^2 - c_g^2) s + \alpha_C^2 s - 2 l(Z) s = 0. $$
Note that because $C$ is effective and because $H$ is ample, $c_h=C.H>0$.

If $c_h^2-c_g^2>0$, then the wall at $t=0$\footnote{It is sufficient to consider the walls at $t=0$ for the following reason: Because the walls for $\cO_S$ in $\Pi_u$ are semicircles \cite{MaciociaWalls}, the wall $\cW(E,\cO_S) \subset \cS_{G,H}$ lies inside the vertical cylinder over $\cW(E,\cO_S)$ at $t=0$.} is an ellipse going through $(0,0)$ and
$$ P_W = \frac{C^2-2l(Z)}{c_g^2-c_h^2} (c_h, c_g), $$
where these are the two points where the tangent line is vertical (see Lemma \ref{PW}).
Therefore, the $s$-value of any point on the ellipse is between $0$ and the $s$-value of $P_W$, which is
$$ \frac{C^2-2l(Z)}{c_g^2-c_h^2} c_h = \frac{c_h^2 - c_g^2 + \alpha_C^2 - 2 l(Z)}{c_g^2 - c_h^2} c_h = -c_h + \frac{\alpha_C^2 - 2l(Z)}{c_g^2-c_h^2}c_h \geq -c_h. $$
This means that the ellipse is contained in the region $s\geq-c_h$.

\begin{figure}[h]
\begin{center}
\includegraphics[scale=.35]{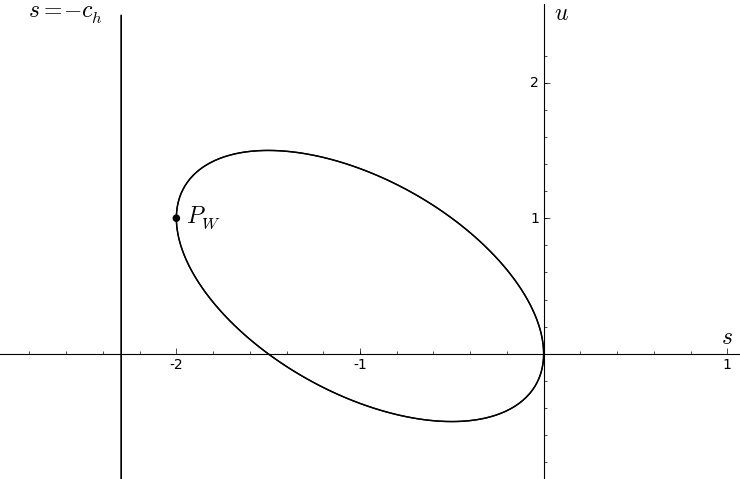}
\label{rank1Ellwall}
\caption{$\cW(\cI_Z(-C),\cO_S)$ at $t=0$ when $c_h^2-c_g^2>0$}
\end{center}
\end{figure}

Because $\cI_Z(-C)\in\cA$, we have that $s<-c_h$, and therefore, $\cI_Z(-C)$ cannot weakly destabilize $\cO_S$.

If $c_h^2-c_g^2=0$, then $C^2\geq0$ implies that $C^2=0$ and $\ch_2(\cI_Z(-C))=-l(Z)\leq0$.
As observed in Section \ref{walls}, where we listed the possible weakly destabilizing walls for $\cO_S$, this wall cannot be a weakly destabilizing wall.
\qed

\begin{remark}
\label{posCnotdestab}
Because $\cO_S$ is Bridgeland stable for $t>>0$ (by Proposition \ref{Olarget}), we have from Proposition \ref{rk1}  that if $C^2\geq0$, then $\beta(\cI_Z(-C))<\beta(\cO_S)$ whenever $\cI_Z(-C)\su\cO_S$ in $\cA$.
\end{remark}

\begin{remark}
\label{negCdestab}
For a curve $C$ of negative self-intersection, we have two possibilities.
\begin{enumerate}

	\item If $C^2<0$ and $c_g^2-c_h^2\leq 0$, then $\cI_Z(-C)$ does not weakly destabilize $\cO_S$ in $\cS_{G,H}$ (note that this is only possible for $S$ with Picard rank $\geq 3$).

	\item If $C^2<0$ and $c_g^2-c_h^2> 0$, then $\cW(\cI_Z(-C),\cO_S)$ is a Left Hyperbola and is a weakly destabilizing wall for $\cO_S$.
\end{enumerate}
\end{remark}


\subsection{Surfaces with no curves of negative self-intersection}

We characterize the stability of $\cO_S$ when $S$ has no curves of negative self-intersection. This implies the first part of Theorem \ref{mainthm}.

\begin{thm}\label{noneg}
If $S$ does not contain any curves of negative self-intersection
and $\sigma_0=(\cZ,\cA_0)\in\stabdivS$ is a stability condition such that $\cO_S\in\cA_0$, then $\cO_S$ is stable with respect to $\sigma_0$.
\end{thm}

\pf
We prove the following statement by induction on the rank of $E$: If $E\su\cO_S$ is a proper subobject in $\cA_0$ for some $\sigma_0=(\cZ,\cA_0)\in\stabdivS$, then $\beta(E)<\beta(\cO_S)$ at $\sigma_0$.

If $E$ has rank $1$, it follows from Lemma \ref{subsofO} that $E$ must be equal to $\cI_Z(-C)$ for some effective curve $C$ and some zero-dimensional scheme $Z$ (with $C$ or $Z$ possibly zero).
We observed in the previous section that these cannot destabilize $\cO_S$ unless $C^2<0$ (see Remark \ref{posCnotdestab}).

Assume now that $E$ has rank $r>1$ and that the result holds for any proper subobject of rank less than $r$ and any stability condition for which the object is indeed a subobject of $\cO_S$.
Choose a pair $G,H$ as above so that $\sigma_0\in\cS_{G,H}$.
From our above study of the walls, we know that the wall $\cW(E,\cO_S)$ is an Ellipse, a Cone, a Left or Right Hyperbola, or a Parabola.
By Remark \ref{In_cA}, we know that $E\su\cO_S\in\cA_0$ iff $\mu_H(\Jl)\leq s_0<\mu_H(\Kl)<0$.
If we had that $\beta(E)\geq\beta(\cO_S)$ at $\sigma_0$, then $E$ would weakly destabilize $\cO_S$.
If the wall was an Ellipse, it would have to intersect $s=\mu_H(\Kl)$ because the wall $\cW(E,\cO_S)$ is connected, and it would pass through a point with $s=s_0<\mu_H(\Kl)$ and $(0,0)$.
If it were a Left Hyperbola, it would have to intersect $s=\mu_H(\Jl)$ because it passes through a point with $s=s_0\geq\mu_H(\Jl)$.
If it were a Right Hyperbola, a Parabola, or a Cone, it would have to intersect $s=\mu_H(\Kl)$ and/or $s=\mu_H(\Jl)$.
Regardless of the type of the wall, we would have that, by Lemma \ref{BertramLemma}, there would exist a proper subobject of $\cO_S$ of higher $\beta$ and lower rank, contradicting our induction hypothesis.
\qed


\subsection{Actual walls are Left Hyperbolas}

There are characteristics of actually destabilizing walls and subobjects that persist regardless of the Picard rank of $S$ or the composition of curves of negative self-intersection within $S$.
The following proposition, which is true in general, will be very useful in our study of surfaces of Picard rank 2 in Section \ref{picrk2}.


\begin{pro}\label{LH}
Let $E\su\cO_S$ in some $\cA$, with quotient $Q$.
If the wall $\cW(E,\cO_S)$ is a destabilizing wall, then it is a Left Hyperbola.
Moreover, $C=c_1(H^0(Q))$ is a curve of negative self-intersection such that the wall $\cW(\cO_S(-C),\cO_S)$ is also a Left Hyperbola, and the wall $\cW(E,\cO_S)$ is inside the wall $\cW(\cO_S(-C),\cO_S)$ for all $|u|>>0$.
\end{pro}

The proof mostly follows from two basic lemmas in which we prove the following results:
\begin{itemize}
\item
Every weakly destabilizing wall that is not a Left Hyperbola is inside a higher (weakly or actually destabilizing) wall that is a Left Hyperbola.
\item
If $E\su\cO_S$ in $\cA$ with quotient $Q$ has a weakly destabilizing wall that is a Left Hyperbola, then either $C=c_1(H^0(Q))$ is a curve of negative self-intersection or the wall is inside a higher (weakly or actually destabilizing) wall that is not a Left Hyperbola.
\end{itemize}

Let $E\su\cO_S$ in $\cA$, and let $Q$ be the quotient.
We use the same notation as in Section \ref{bertram} for the Harder-Narasimhan filtrations of $E$ and $H^{-1}(Q)$ with respect to the Mumford slope $\mu_H$.

\begin{lem}\label{notLH}
Let $E\su\cO_S$ in some $\cA$, and suppose that the wall $\cW(E,\cO_S)$ is a weakly destabilizing wall that is not a Left Hyperbola.
Then, for some $E_i$ in the Harder-Narasimhan filtration of $E$, the wall $\cW(E_i,\cO_S)$ is a Left Hyperbola such that the following is true:
If there exists a stability condition $\sigma_0=\sigma_{s_0H+u_0G,t_0H}=(\cZ,\cA_0)$ such that $E\su\cO_S$ in $\cA_0$, then $E_i\su\cO_S$ in $\cA_0$, and the wall $\cW(E,\cO_S)$ is inside the wall $\cW(E_i,\cO_S)$ at $u_0$.
In particular, $E$ cannot actually destabilize $\cO_S$ anywhere.
\end{lem}

\pf
We prove this by induction on the number of terms $n$ in the Harder-Narasimhan filtration of $E$.

Suppose that $n=1$, i.e., that $E$ is $\mu_H$-semistable.
Because the wall $\cW(E,\cO_S)$ is a weakly destabilizing wall, there exists a stability condition $\sigma_0=(Z,\cA_0)$ such that $\beta(E)\geq\beta(\cO_S)$ at $\sigma_0$.
By Remark \ref{In_cA}, this implies that there exist a point on the wall with $s$ such that $s=s_0<\mu_H(E)<0$.
Therefore, the wall $\cW(E,\cO_S)$ cannot be a weakly destabilizing wall without being a Left Hyperbola because all other types of walls would not be able to contain this point with $s$-value $s_0<\mu_H(E)$ and the point $(0,0)$ without intersecting the line $s=\mu_H(E)$, which is not possible \cite{MaciociaWalls}.

Now, let $n>1$, and assume that the statement is true for all subobjects of $\cO_S$ that have a weakly destabilizing wall that is not a Left Hyperbola and have a Harder-Narasimhan filtration of length $<n$.

Because the wall $\cW(E,\cO_S)$ is not a Left Hyperbola, it will intersect the line $s=\mu_H(\Kl)$ by the same argument as that of the $n=1$ case.
Consider a value of $u$ such that $\cW(E,\cO_S)\cap\Pi_u$ intersects $s=\mu_H(\Kl)$ at $t>0$.
By Lemma \ref{BertramLemma}, we have that $\beta(E_{n-1})>\beta(E)$ at the stability conditions $\sigma_{s,u}$ on $\cW(E,\cO_S)\cap\Pi_u$ such that $s<\mu_H(\Kl)$.
Therefore, at those stability conditions, we have $\beta(E_{n-1})>\beta(E)=\beta(\cO_S)$, and the wall $\cW(E_{n-1},\cO_S)$ is also a weakly destabilizing wall.

If the wall $\cW(E_{n-1},\cO_S)$ is not a Left Hyperbola, we can conclude by induction that there exists an $E_i$ such that the wall $\cW(E_i,\cO_S)$ is a Left Hyperbola and that the wall $\cW(E_{n-1},\cO_S)$ is inside the wall $\cW(E_i,\cO_S)$ for all $u_0$ for which there exists a stability condition $\sigma_0=\sigma_{s_0H+u_0G,t_0H}=(\cZ,\cA_0)$ with $E_{n-1}\su\cO_S$ in $\cA_0$.
If the wall $\cW(E_{n-1},\cO_S)$ is a Left Hyperbola, let $i=n-1$.

Now, let $\sigma_0=\sigma_{s_0H+u_0G,t_0H}=(\cZ,\cA_0)$ be a stability condition such that $E\su\cO_S$ in $\cA_0$.
We need to prove that the wall $\cW(E,\cO_S)$ is inside the wall $\cW(E_i,\cO_S)$ at $u_0$.
First, notice that if $E\su\cO_S$ in $\cA_0$, then $E_i\su E\su\cO_S$ in $\cA_0$ because if $E\in\cA_0$, then $E_i\su E\in\cA_0$.

Now, suppose that $E$ weakly destabilizes $\cO_S$ at $\sigma_0$ (i.e. $\beta(E)\geq\beta(\cO_S)$ at $\sigma_0$), and suppose moreover that $u_0>0$ (the proof for $u_0<0$ is similar).
Because the wall $\cW(E,\cO_S)$ is not a Left Hyperbola, its intersection with the planes $\Pi_u$ is non-empty for all $0\leq u\leq u_0$.
On the other hand, because the wall $\cW(E_i,\cO_S)$ is a Left Hyperbola, its region in the $s<0$ half plane does not reach the $s$-axis.
Therefore, for small positive values of $u$, $\cW(E_i,\cO_S)\cap\Pi_u$ is empty in the $s<0$ region, and $\cW(E,\cO_S)\cap\Pi_u$ is non-empty.

By Lemma \ref{BertramLemma} and the induction hypothesis, we have that
$$ \cW(E,\cO_S) \cap \Pi_u \preceq \cW(E_{n-1},\cO_S) \cap \Pi_u \preceq \cW(E_i,\cO_S) \cap \Pi_u $$
for every value of $u$ such that $\cW(E,\cO_S)\cap\Pi_u$ intersects $s=\mu_H(\Kl)$ at $t>0$.
The statement remains true at $u_1$ by continuity, where $u_1$ is the smallest value of $u$ such that $\cW(E,\cO_S)$ intersects $s=\mu_H(\Kl)$ (the intersection will be at $t=0$).

Because the wall $\cW(E_i,\cO_S)$ is inside the wall $\cW(E,\cO_S)$ for small positive values of $u$, and the wall $\cW(E,\cO_S)$ is inside the wall $\cW(E_i,\cO_S)$ at $u_1$, there must exist a value $u_2$ with $0<u_2<u_1$ such that
$ \cW(E,\cO_S)\cap\Pi_{u_2} = \cW(E_i,\cO_S)\cap\Pi_{u_2} $
(recall that, for each $u$, the intersections of the walls with the plane $\Pi_u$ are nested semicircles).
Moreover,
$$ \cW(E,\cO_S) \cap \Pi_u \preceq \cW(E_i,\cO_S) \cap \Pi_u $$
for all $u\geq u_2$ by Lemma \ref{geometry}.
Because $E\in\cA_0$, we have that $s_0<\mu_H(\Kl)$.
However, the wall $\cW(E,\cO_S)$ is not a Left Hyperbola; therefore, $u_0\geq u_1>u_2$, and thus, the wall $\cW(E,\cO_S)$ is inside the wall $\cW(E_i,\cO_S)$ at $u_0$.
\qed

\begin{figure}[h]
\begin{center}
\includegraphics[scale=.35]{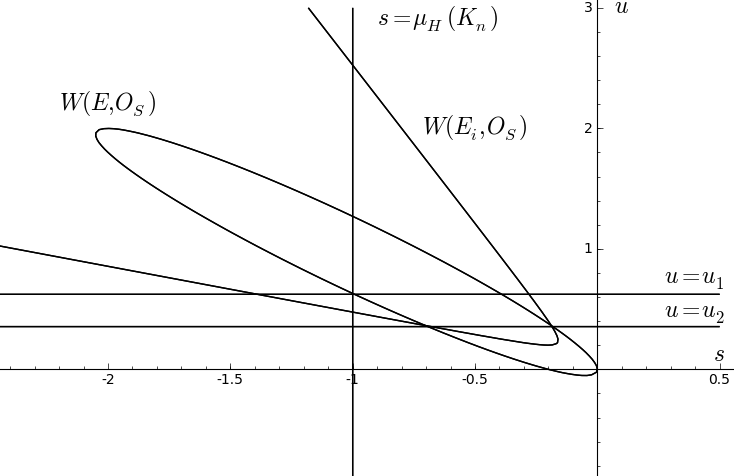}
\caption{The case of an ellipse in the proof of Lemma \ref{notLH}}
\end{center}
\end{figure}

\begin{lem}\label{LemmaLH}
Let $E\su\cO_S$ in some $\cA$ with quotient $Q$, and suppose that the wall $\cW(E,\cO_S)$ is a weakly destabilizing Left Hyperbola.
Then,
\begin{itemize}
\item
Either $C=c_1(H^0(Q))$ is a curve of negative self-intersection and the wall $\cW(\cO_S(-C),\cO_S)$ is a Left Hyperbola, or
\item
For some $F_j$ in the Harder-Narasimhan filtration of $H^{-1}(Q)$, the wall $\cW(E/F_j,\cO_S)$ is not a Left Hyperbola, and the following is true:
If there exists a stability condition $\sigma_0=\sigma_{s_0H+u_0G,t_0H}=(\cZ,\cA_0)$ such that $E\su\cO_S$ in $\cA_0$, then $E/F_j\su\cO_S$ in $\cA_0$, and the wall $\cW(E,\cO_S)$ is inside the wall $\cW(E/F_j,\cO_S)$ at $u_0$.
In particular, $E$ cannot actually destabilize $\cO_S$ anywhere.
\end{itemize}
\end{lem}

\pf
We prove this by induction on the number of terms $m$ in the Harder-Narasimhan filtration of $H^{-1}(Q)$ (including the case $m=0$ corresponding to $H^{-1}(Q)=0$).

If $m=0$ and $H^{-1}(Q)=0$, then $E=\cI_Z(-C)$, and by Proposition \ref{rk1}, $C=c_1(H^0(Q))$ must be a curve of negative self-intersection for the wall $\cW(E,\cO_S)$ to be a weakly destabilizing wall.
Moreover, the wall must be a Left Hyperbola by Remark \ref{negCdestab}.

Now, let $m>0$, and assume that the statement is true for all subobjects $E'$ of $\cO_S$ that have a weakly destabilizing wall that is a Left Hyperbola and a Harder-Narasimhan filtration for $H^{-1}(Q')$ of length $<m$, where $Q'$ is the quotient of $\cO_S$ by $E'$.

Because the wall $\cW(E,\cO_S)$ is a Left Hyperbola, it will intersect the line $s=\mu_H(\Jl)=\mu_H(F_1)$.
Consider a value of $u$ such that $\cW(E,\cO_S)\cap\Pi_u$ intersects $s=\mu_H(F_1)$ at $t>0$.
By Lemma \ref{BertramLemma}, $\beta(E/F_1)>\beta(E)$ at those stability conditions.
Therefore, the wall $\cW(E/F_1,\cO_S)$ is also a weakly destabilizing wall.

Let $Q_1$ be the quotient for $E/F_1\su\cO_S$.
By Remark \ref{emodjremark}, we know that $H^0(Q_1)=H^0(Q)$ and $H^{-1}(Q_1)=H^{-1}(Q)/F_1$.
In particular, the Harder-Narasimhan filtration of $H^{-1}(Q_1)$ is simply $F_j/F_1$ ($1\leq j\leq m$), which has length $m-1$.
Moreover, the quotients $(E/F_1)/(F_j/F_1)$ are the quotients $E/F_j$.

If the wall $\cW(E/F_1,\cO_S)$ is a Left Hyperbola, we have by induction that at least one of the following two options is true:

(1) $C=c_1(H^0(Q_1))$ is a curve of negative self-intersection, and the wall $\cW(\cO_S(-C),\cO_S)$ is a Left Hyperbola.
Because $H^0(Q_1)=H^0(Q)$, we are done.

(2) There exists an $F_j/F_1$ in the Harder-Narasimhan filtration of $H^{-1}(Q_1)$ such that the wall $\cW(E/F_j,\cO_S)$ is not a Left Hyperbola, and the wall $\cW(E/F_1,\cO_S)$ is inside the wall $\cW(E/F_j,\cO_S)$ for all $u_0$ for which there exists a stability condition $\sigma_0=\sigma_{s_0H+u_0G,t_0H}=(\cZ,\cA_0)$ with
$E/F_1\su\cO_S$ in $\cA_0$, with $E/F_j\su\cO_S$ there.

If the wall $\cW(E/F_1,\cO_S)$ is not a Left Hyperbola, let $j=1$.

Now, let $\sigma_0=\sigma_{s_0H+u_0G,t_0H}=(\cZ,\cA_0)$ be a stability condition such that $E\su\cO_S$ in $\cA_0$.
We need to prove that the wall $\cW(E,\cO_S)$ is inside the wall $\cW(E/F_j,\cO_S)$ at $u_0$.

First, notice that if $E\su\cO_S$ in $\cA_0$, then $E/F_j\su\cO_S\in\cA_0$.

As in the previous proof, assume that $E$ weakly destabilizes $\cO_S$ at $\sigma_0$ and that $u_0>0$.
We then know that $\cW(E/F_j,\cO_S)\cap\Pi_u$ is not empty for small positive values of $u$, and $\cW(E,\cO_S)\cap\Pi_u$ is empty in the region $s<0$ for $u$ positive and sufficiently small.

Let $u_1$ be the largest value of $u$ such that the wall $\cW(E,\cO_S)$ intersects $s=\mu_H(\Jl)$ (the intersection will be at $t=0$).
As in the previous proof, we can use Lemma \ref{BertramLemma}, the induction hypothesis, and continuity to show that
$$ \cW(E,\cO_S) \cap \Pi_{u_1} \preceq \cW(E/F_1,\cO_S) \cap \Pi_{u_1} \preceq \cW(E/F_j,\cO_S) \cap \Pi_{u_1}. $$

Because $\cW(E,\cO_S)\cap\Pi_u$ is inside $\cW(E/F_j,\cO_S)\cap\Pi_u$ for small positive values of $u$ and at $u_1$, it must be inside of it for all $0<u\leq u_1$.
Otherwise, there would have to exist two values of $u$ between $0$ and $u_1$ whereby
$ \cW(E,\cO_S)\cap\Pi_u = \cW(E/F_j,\cO_S)\cap\Pi_u, $
which is not possible by Lemma \ref{geometry}.

Because $E\su\cO_S\in\cA_0$, we have that $s_0\geq\mu_H(\Jl)$.
Moreover, $u_0\leq u_1$ because $\cW(E,\cO_S)$ is a Left Hyperbola, and therefore, the wall $\cW(E,\cO_S)$ is inside the wall $\cW(E/F_j,\cO_S)$ at $u_0$.
\qed

\begin{figure}[h]
\begin{center}
\includegraphics[scale=.35]{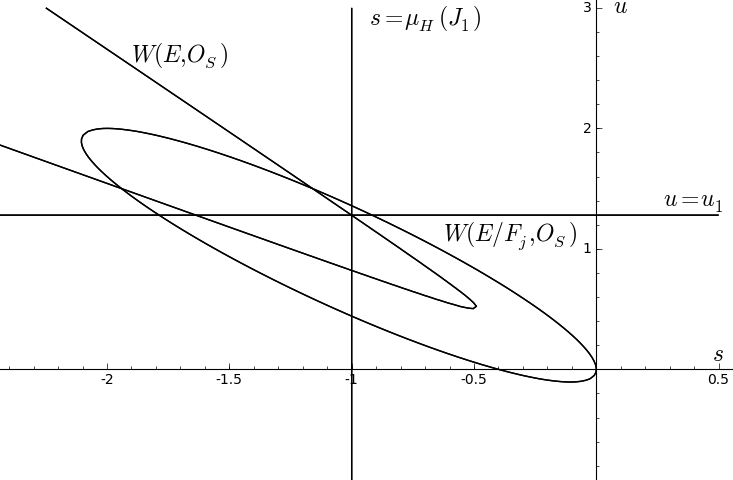}
\caption{A possible case for the proof of Lemma \ref{LemmaLH}}
\end{center}
\end{figure}

We can now prove Proposition \ref{LH}.
The only part of the statement of the proposition that does not directly follow from the two lemmas is the part where we claim that the wall $\cW(E,\cO_S)$ is inside the wall $\cW(\cO_S(-C),\cO_S)$ for all $|u|>>0$.

We prove something a bit stronger:

\begin{lem}\label{LHuLarge}
Let $E\su\cO_S$ in $\cA$ be an actually destabilizing object with quotient $Q$.
For all $F_j$ in the Harder-Narasimhan filtration of $H^{-1}(Q)$, the wall $\cW(E/F_j,\cO_S)$ is a Left Hyperbola, and the wall $\cW(E,\cO_S)$ is inside the wall $\cW(E/F_j,\cO_S)$ for all $|u|>>0$.
\end{lem}

The statement of the proposition follows from this lemma because $F_m=H^{-1}(Q)$ and $E/F_m=\cI_Z(-C)$ for some zero-dimensional scheme $Z$.
Then, for $|u|>>0$, the wall $\cW(E,\cO_S)$ must be inside the wall $\cW(\cI_Z(-C),\cO_S)$.
The result then follows from the fact that the wall $\cW(\cI_Z(-C),\cO_S)$ is always inside the wall $\cW(\cO_S(-C),\cO_S)$.
Indeed, if $\cI_Z(-C)\su\cO_S(-C)$ in $\cA$, then $\cZ(\cI_Z(-C))=\cZ(\cO_S(-C))+l(Z)$, where $l(Z)$ is the length of $Z$, and $\beta(\cI_Z(-C))<\beta(\cO_S(-C))$.
Another way to see this is through the fact that the quotient of $\cO_S(-C)$ by $\cI_Z(-C)$ is a torsion sheaf $T_Z$ supported on $Z$, which satisfies $\cZ(T_Z)=-l(Z)$ and has maximal phase $\beta(T_Z)=\infty$.
This means that $\cI_Z(-C)$ can never destabilize $\cO_S(-C)$, and the wall $\cW(\cI_Z(-C),\cO_S)$ is always inside the wall $\cW(\cO_S(-C),\cO_S)$ (otherwise, we would have $\beta(\cI_Z(-C))>\beta(\cO_S)>\beta(\cO_S(-C))$ in the region between the two walls).

\vspace{.1in}

{\it Proof} (of Lemma \ref{LHuLarge}).
Assume that $E$ actually destabilizes $\cO_S$ at a stability condition $\sigma_0=\sigma_{s_0H+u_0G,t_0H}=(\cZ,\cA_0)$ and that $u_0>0$ (the proof for $u_0<0$ is very similar).
We know that all of the walls $\cW(E/F_j,\cO_S)$ must be Left Hyperbolas by the proof of the previous lemma.
Indeed, if they were not Left Hyperbolas, then there would  exist a $j$ such that the wall $\cW(E,\cO_S)$ would be inside the wall $\cW(E/F_j,\cO_S)$ at $u_0$, making it impossible for $E$ to actually destabilize $\cO_S$ there.

Let $u_1$ be the largest value of $u$ such that $\cW(E,\cO_S)$ intersects $s=\mu_H(\Jl)$ (the intersection will be at $t=0$).
Then, we know that the region where $E\su\cO_S$ in $\cA$ and $\beta(E)>\beta(\cO_S)$ is contained within $s\geq\mu_H(F_1)$ and $0<u\leq u_1$.
Moreover, as above, we have that
$ \cW(E,\cO_S) \cap \Pi_{u_1} \preceq \cW(E/F_1,\cO_S) \cap \Pi_{u_1}. $

For $E$ to actually destabilize $\cO_S$ at $\sigma_0$, we must have that $0<u_0<u_1$, and
$ \cW(E/F_1,\cO_S) \cap \Pi_{u_0} \preceq \cW(E,\cO_S) \cap \Pi_{u_0}. $
Using Lemma \ref{geometry} as above, we see that because the wall $\cW(E/F_1,\cO_S)$ is inside the wall $\cW(E,\cO_S)$ at $u_0$, and because the wall $\cW(E,\cO_S)$ is inside the wall $\cW(E/F_1,\cO_S)$ at $u_1>u_0$, the wall $\cW(E,\cO_S)$ must remain inside the wall $\cW(E/F_1,\cO_S)$ for all $u\geq u_1$.
Thus, if $m=1$, we are done.

If $m>1$, let $u_2>u_1$ be the largest value of $u$ such that $\cW(E/F_1,\cO_S)$ intersects $s=\mu_H(F_2/F_1)$.
Then,
$ \cW(E,\cO_S) \cap \Pi_{u_2} \preceq \cW(E/F_1,\cO_S) \cap \Pi_{u_2} \preceq \cW(E/F_2,\cO_S) \cap \Pi_{u_2} $
as above.

Because the wall $\cW(E/F_2,\cO_S)$ is inside the wall $\cW(E,\cO_S)$ at $u_0$, and because the wall $\cW(E,\cO_S)$ is inside the wall $\cW(E/F_1,\cO_S)$ at $u_2>u_1>u_0$, we know that the wall $\cW(E,\cO_S)$ must remain inside for all $u\geq u_2$.

Continuing in a similar manner proves the statement for all $1\leq j\leq m$.
\qed


\subsection{Surfaces of Picard rank 2}
\label{picrk2}

We now restrict our attention to surfaces of Picard rank 2, whereby we can more precisely describe the actually destabilizing walls for $\cO_S$.
We pause to illuminate a fact that is helpful in this situation.

\begin{remark}
\label{stayinSHG}
Let $S$ have Picard rank 2.
For any line bundle $\cO_S(D')$ and $\sigma_{D,H}\in\cS_{G,H}$, we have $\cO_S(D')\otimes\sigma_{D,H}=\sigma_{D+D',H}\in\cS_{G,H}$.
\end{remark}

\pf
This is essentially Lemma \ref{reductionlemma}.
The important fact to note in the case of Picard rank 2, however, is that if $\sigma_{D,H}\in\cS_{G,H}$, then $\sigma_{D+D',H}$ is also in $\cS_{G,H}$ for all $D'$.
This would not necessarily be true for a surface of Picard rank $>2$ because the divisor $D'$ may not be a linear combination of $G$ and $H$.
\qed

Let $S$ be a surface of Picard rank 2, and let $H$ and $G$ be as above.
Moreover, let $C_1$ and $C_2$ be the generators of the cone of effective curves on $S$.
Because $H$ is ample, we must have that $H=eC_1+fC_2$ for some $e,f>0$.
Then, $H.G=0$ implies that $fC_2.G=-eC_1.G$.
Therefore, $(C_1.G)\cdot(C_2.G)<0$.
Assume that $C_1.G>0$ and $C_2.G<0$.

We observed in Proposition \ref{LH} that every destabilizing wall is a Left Hyperbola that has to be inside a rank-1 weakly destabilizing wall for $|u|>>0$.
In the case of Picard rank 2, we can make $|u|>>0$ more precise as follows:

\begin{pro}\label{conjevidencePicRk2}
If $u\geq C_1.G$, then $\cO_S$ is only destabilized by $\cO_S(-C_1)$, and if $u\leq C_2.G$, then $\cO_S$ is only destabilized by $\cO_S(-C_2)$.
\end{pro}

\pf
We prove the following statement by induction on the rank of $E$:
If $E\su\cO_S$ for some stability condition $\sigma_0=\sigma_{s_0H+u_0G,t_0H}=(\cZ,\cA_0)$ with $u_0\geq C_1.G$, and $\beta(E)\geq\beta(\cO_S)$ at $\sigma_0$, then the wall $\cW(E,\cO_S)$ is inside the wall $\cW(\cO_S(-C_1),\cO_S)$ at $u_0$.
(A similar proof would work with $C_2$ in place of $C_1$ if $u_0\leq C_2.G$.)

If the rank of $E$ is 1 and $\beta(E)\geq\beta(\cO_S)$ at $\sigma_0$, then $E=\cI_Z(-C)$ for some curve $C$ of negative self-intersection and some zero-dimensional scheme $Z$.
We observed above, between the statement and the proof of Lemma \ref{LHuLarge}, that $\beta(\cO_S(-C))\geq\beta(\cI_Z(-C))$.
Therefore, we have that $\beta(\cO_S(-C))\geq\beta(\cI_Z(-C))\geq\beta(\cO_S)$ at $\sigma_0$ with $u_0\geq C_1.G>0$, and we must have that $C=aC_1+bC_2$ with $a>0$.
This gives us $\cO_S(-C)\su\cO_S(-C_1)$.

Note that the wall $\cW(\cO_S(-C),\cO_S)$ is inside the wall $\cW(\cO_S(-C_1),\cO_S)$ for $u>>0$ because of the slopes of their asymptotes (see Section \ref{walls}), and the former is inside of the latter at $u=C_1.G$ because $\cO_S(-C_1)$ is always Bridgeland-stable there.
Therefore, the former must be inside of the latter for all $u\geq C_1.G$ by Lemma \ref{geometry}.

Assume now that $E$ has rank $r>1$ and that the statement is true for all subobjects of $\cO_S$ of rank $<r$.

By Lemma \ref{notLH}, we can assume that the wall $\cW(E,\cO_S)$ is a Left Hyperbola.
Moreover, if we let $C=c_1(H^0(Q))$, where $Q$ is the quotient of $\cO_S$ by $E$ in $\cA$, we have that the wall $\cW(\cO_S(-C),\cO_S)$ is also a Left Hyperbola and that the wall $\cW(E,\cO_S)$ is inside the wall $\cW(\cO_S(-C),\cO_S)$ for $u>>0$.
From the rank 1 case, we know that the wall $\cW(\cO_S(-C),\cO_S)$ is inside the wall $\cW(\cO_S(-C_1),\cO_S)$ for all $u\geq C_1.G$.

If the wall $\cW(E,\cO_S)$ is inside the wall $\cW(\cO_S(-C_1),\cO_S)$ at $u_0$, we are done.
Suppose that it is not.
Then, we have by Lemma \ref{geometry} that the wall $\cW(\cO_S(-C_1),\cO_S)$ is inside the wall $\cW(E,\cO_S)$ for all $u\leq u_0$. Denote this statement by ($\star$).

If $\cW(E,\cO_S)\cap\Pi_{u_0}$ intersects $s=\mu_H(\Kl)$ (using our usual notation for the Harder-Narasimhan filtration from Section \ref{bertram}), then
$$ \cW(E,\cO_S) \cap \Pi_{u_0} \preceq
\cW(E_{n-1},\cO_S) \cap \Pi_{u_0} \preceq \cW(\cO_S(-C_1),\cO_S) \cap \Pi_{u_0} $$
by Lemma \ref{BertramLemma} and the induction hypothesis, contradicting ($\star$).
Therefore, $\cW(E,\cO_S)\cap\Pi_{u_0}$ does not intersect $s=\mu_H(\Kl)$, and because $u_0\geq C_1.G$, ($\star$) implies that the wall $\cW(\cO_S(-C_1),\cO_S)$ is inside the wall $\cW(E,\cO_S)$ at $u=C_1.G$.
Because $\cO_S(-C_1)$ is Bridgeland stable at $u=C_1.G$, this can only occur if $E$ is not a subobject of $\cO_S(-C_1)$ in $\cA$ at the points on the wall $\cW(E,\cO_S)\cap\Pi_{C_1.G}$.
This means that the wall $\cW(E,\cO_S)\cap\Pi_u$ had to intersect $s=\mu_H(\Kl)$ for some $C_1.G\leq u<u_0$.

Let $u_1$ be a value of $u$ with $C_1.G\leq u_1<u_0$ where the wall $\cW(E,\cO_S)$ intersects $s=\mu_H(\Kl)$.
Then,
$ \cW(E,\cO_S) \cap \Pi_{u_1} \preceq
\cW(E_{n-1},\cO_S) \cap \Pi_{u_1} \preceq \cW(\cO_S(-C_1),\cO_S) \cap \Pi_{u_1} $
by Lemma \ref{BertramLemma} and the induction hypothesis, contradicting ($\star$).
Thus, the wall $\cW(E,\cO_S)$ is inside the wall $\cW(\cO_S(-C_1),\cO_S)$ at $u_0$.
\qed




\vspace{.1in}

We now prove the second part of our main Theorem \ref{mainthm}, i.e., that Conjecture \ref{mainconj} is true for surfaces of Picard Rank 2 that only have one irreducible curve of negative self-intersection. We first give a lemma and then prove a result that is stronger than the conjecture in this situation.


\begin{lem}\label{onenegcurvelem}
Let $S$ be a surface of Picard Rank $2$.
Assume that the cone of effective curves is generated by $C_1$ and $C_2$, that $C_1.G>0$, and that $C_1$ is the only irreducible curve in $S$ of negative self-intersection.
Let $C$ be an effective curve.
If $C.G<0$, then $C^2\geq0$.
\end{lem}

\pf
We observed above that we can write $H=eC_1+fC_2$ with $e,f>0$, and obtain $fC_2.G=-eC_1.G$.
Because $(eC_1+fC_2)^2=H^2>0$, we have that $2efC_1.C_2>-e^2C_1^2-f^2C_2^2$.

Let $C$ be an effective curve such that $C.G<0$.
Then, $C=aC_1+bC_2$ with $a,b\geq0$, and $0>fC.G=afC_1.G+bfC_2.G=(af-be)C_1.G$.
Therefore, $af-be<0$ because $C_1.G>0$.

This gives us that
$ efC^2 = a^2ef C_1^2 + b^2ef C_2^2 + 2abef C_1.C_2 > a^2ef C_1^2 + b^2ef C_2^2 - abe^2 C_1^2 - abf^2 C_2^2 = ae (af - be) C_1^2 + bf (be - af) C_2^2 \geq 0, $
and $C^2\geq0$.
\qed

\begin{pro}\label{onenegcurve}
Let $S$ be a surface of Picard Rank $2$.
Assume that the cone of effective curves is generated by $C_1$ and $C_2$, that $C_1.G>0$ and that $C_1$ is the only irreducible curve in $S$ of negative self-intersection.
Then, $\cO_S$ is only destabilized by $\cO_S(-C_1)$.
\end{pro}

\pf
Let us first note that because $C_1$ is the only irreducible curve of negative self-intersection, we must have $C_2^2\geq0$.
Therefore, $\cO_S(-C_2)$ does not weakly destabilize $\cO_S$, and the wall $\cW(\cO_S(-C_2),\cO_S)$ must be empty in the region where $\cO_S(-C_2)\in\cA$.

We prove the following statement by induction on the rank of $E$:
If $E\su\cO_S$ for some stability condition $\sigma_0=\sigma_{s_0H+u_0G,t_0H}=(\cZ,\cA_0)$, and $\beta(E)\geq\beta(\cO_S)$ at $\sigma_0$, then the wall $\cW(E,\cO_S)$ is inside the wall $\cW(\cO_S(-C_1),\cO_S)$ at $u_0$.

If $u_0\geq C_1.G$, we already know this statement to be true by the proof of Proposition \ref{conjevidencePicRk2}.
Assume now that $u_0<C_1.G$.

If the rank of $E$ is 1, and if $\beta(E)\geq\beta(\cO_S)$ at $\sigma_0$, then $E=\cI_Z(-C)$ for some curve $C$ of negative self-intersection and some zero-dimensional scheme $Z$.
Suppose, by contradiction, that the wall $\cW(E,\cO_S)$ is not inside the wall $\cW(\cO_S(-C_1),\cO_S)$ at $u_0$.
Choose a stability condition $\sigma'$ on the semicircle $\cW(E,\cO)\cap\Pi_{u_0}$.
Then, at $\sigma'$, we have $\beta(\cO_S(-C))\geq\beta(\cI_Z(-C))>\beta(\cO_S(-C_1))$.
Therefore, at $\cO_S(C_1)\otimes\sigma'$, we have $\beta(\cO_S(-C+C_1))>\beta(\cO_S)$ (see the proof of Lemma \ref{reductionlemma}).
This means that $\cO_S(-C+C_1)$ weakly destabilizes $\cO_S$ at $\cO_S(C_1)\otimes\sigma'$.
By Proposition \ref{rk1}, we must have that $(C-C_1)^2<0$, and by Lemma \ref{onenegcurvelem}, $(C-C_1).G\geq0$.
This implies that the wall $\cW(\cO_S(-C+C_1),\cO_S)$ is a Left Hyperbola that can only weakly destabilize $\cO_S$ in the region $u>0$, but we have that the value of $u$ at $\cO_S(C_1)\otimes\sigma'$ is equal to $u_0-C_1.G<0$.
This proves that the statement is true if $E$ has rank $1$.

Assume now that $E$ has rank $r>1$ and that the statement is true for all subobjects of $\cO_S$ of rank $<r$.

As in the proof of Proposition \ref{conjevidencePicRk2}, we let $C=c_1(H^0(Q))$, and we can then assume that $\cW(E,\cO_S)$ and $\cW(\cO_S(-C),\cO_S)$ are Left Hyperbolas and that the wall $\cW(E,\cO_S)$ is inside the wall $\cW(\cO_S(-C),\cO_S)$ for $u>>0$.
In addition, by the rank 1 case, we know that the wall $\cW(\cO_S(-C),\cO_S)$ is inside the wall $\cW(\cO_S(-C_1),\cO_S)$ for all $u$.

If the wall $\cW(E,\cO_S)$ is inside the wall $\cW(\cO_S(-C_1),\cO_S)$ at $u_0$, we are done.
Suppose that this was not the case.
Using our usual notation for the Harder-Narasimhan filtration, we know from the proof of Proposition \ref{conjevidencePicRk2} that $\cW(E,\cO_S)\cap\Pi_{u_0}$ cannot intersect $s=\mu_H(\Kl)$.
Choose a $\sigma'$ in the semicircle $\cW(E,\cO_S)\cap\Pi_{u_0}$.
Then, $\beta(E)>\beta(\cO_S(-C_1))$ at $\sigma'$ implies that $\beta(E(C_1))>\beta(\cO_S)$ at $\cO_S(C_1)\otimes\sigma'$.
Therefore, $E(C_1)$ weakly destabilizes $\cO_S$ at $\cO_S(C_1)\otimes\sigma'$.

Let $u_1$ be the $u$ value of $\cO_S(C_1)\otimes\sigma'$, i.e., let $u_1=u_0-C_1.G<0$.
Consider the highest semi-circular wall at $u_1$, corresponding to an object $E'\su\cO_S$.
Because $E'$ actually destabilizes $\cO_S$, we have that the wall $\cW(E',\cO_S)$ must be a Left Hyperbola.
However, the proof of Proposition \ref{conjevidencePicRk2} shows that the wall $\cW(E',\cO_S)$ would have to be inside the wall $\cW(\cO_S(-C_2),\cO_S)$ for all $u\leq C_2.G$.
However, this cannot be true because the wall $\cW(E,\cO_S)$ is a Left Hyperbola, and the wall $\cW(\cO_S(-C_2),\cO_S)$ is empty.
\qed


\section{Bridgeland stability of $\cO_S[1]$}
\label{Oshiftedby1}

We now proceed to the study of the stability of $\cO_S[1]$.
This can be addressed via duality, except for the stability conditions $\sigma_{D,H}$ with $D=uG$, i.e., $s=0$. 

Note that, given a Bridgeland stability condition $\sigma$, $\cO_S[1]\in\cA$ iff $s\geq0$.

\subsection{Subobjects of $\cO_S[1]$}

Let $\sigma$ be a Bridgeland stability condition, and let $E$ be a proper subobject of $\cO_S[1]$ in $\cA$.
We have a short exact sequence $0\fun E\fun\cO_S[1]\fun Q'\fun0$ in $\cA$ for some $Q'\in\cA$.
The long exact sequence in cohomology is
$$ 0 \lfun H^{-1}(E) \lfun \cO_S \lfun H^{-1}(Q') \lfun H^0(E) \lfun 0, $$
and therefore, $Q'=H^{-1}(Q')[1]$ is the shift of a sheaf.
In addition, because $H^{-1}(Q')\in\cF$ is torsion free, we have that either $H^{-1}(E)=\cO_S$ or $H^{-1}(E)=0$.
However, if $H^{-1}(E)=\cO_S$, then $H^{-1}(Q')=H^0(E)$, and this is not possible because the first sheaf is in $\cF$ and the second one is in $\cT$.
Therefore, $H^{-1}(E)=0$, and $E=H^0(E)$ is a sheaf.
We will denote $H^{-1}(Q')$ by $Q$.

To summarize, if $E\su\cO_S[1]$ is a proper subobject in $\cA$, then $E$ is a sheaf in $\cT$, and the quotient is of the form $Q[1]$ for some sheaf $Q\in\cF$.
We have a short exact sequence of sheaves
$$ 0 \lfun \cO_S \lfun Q \lfun E \lfun 0. $$

\subsection{The $s=0$ case}

Let us start by proving that $\cO_S[1]$ is Bridgeland stable when $s=0$.

\begin{lem}
\label{sequalszero}
If $s=0$, $\cO[1]$ has no proper subobjects in $\cA$ and is therefore Bridgeland stable.
\end{lem}

\pf
Let $E\su\cO_S[1]$ be a proper subobject of $\cO_S[1]$ in $\cA$ with quotient $Q[1]$ as above.
Because $s=0$, we have that $\Im \cZ(\cO_S[1])=0$.
Then, $\Im \cZ(\cO_S[1])=\Im \cZ(E)+\Im \cZ(Q[1])$, and they all have non-negative imaginary parts; therefore, we must have that $\Im \cZ(E)=\Im \cZ(Q[1])=0$.
We claim that $E$ must be a torsion sheaf supported in dimension $0$.
Recall that if $\ch(E)=(r,d_hH+d_gG+\alpha,c)$, then $\Im \cZ(E)=td_h-rst$.
Because $s=0$, we must have in this case that $td_h=0$.
If $E$ has positive rank, then $E\in\cT$ would imply that $d_h>0$, which is not possible.
If $E$ has rank $0$, then $d_h=0$ would imply that it has to be supported in dimension $0$ because $H$ is ample.
Therefore, we have a short exact sequence of sheaves $0\fun\cO_S\fun Q\fun E\fun0$, with $E$ a torsion sheaf supported in dimension $0$.
However, this cannot occur unless the sequence splits, in which case $Q$ would have torsion, which is impossible.
This proves that, if $s=0$, $\cO_S[1]$ cannot have proper subobjects and is Bridgeland stable.
\qed

\vspace{.1in}

Herein, we assume that $s>0$.

\subsection{Duality}

The following duality result allows us to apply the results concerning the stability of $\cO_S$ to that of $\cO_S[1]$. The result follows as in \cite[Lemma 3.2]{Mart} with a slightly different choice of functor. Specifically, we consider the functor $E \mapsto E^\vee := R\mathcal{H}om(E,\cO_S)[1]$. Note that $\cO_S^\vee = \cO_S[1]$ and vice versa.

\begin{lem}
Let $D$ be a divisor with $D.H<0$. Then, $\cO_S$ is $\sigma_{D,H}$-(semi)stable if and only if $\cO_S[1]$ is $\sigma_{-D,H}$-(semi)stable.
\end{lem}

\pf
This follows from \cite[Lemma 3.2(d)]{Mart}.
\qed

\vspace{.1in}

Note that if $\cO_S(-C)\su\cO_S$ destabilizes $\cO_S$ at $\sigma_{D,H}$, then applying $(\_)^\vee$ shows that the quotient $\cO_S[1] \sur \cO_S(C)[1]$ destabilizes $\cO_S[1]$ at $\sigma_{-D,H}$.
The kernel of the map $\cO_S[1] \sur \cO_S(C)[1]$ in $\cA$ is $\cO_S(C)|_C$.
Thus, Theorem \ref{noneg}, Proposition \ref{onenegcurve}, and Lemma \ref{sequalszero} yield the following result.

\begin{pro}
\label{O1stab}
\mbox{}
\begin{enumerate}
\item
If $S$ does not contain any curves of negative self-intersection, and if $\sigma=(\cZ,\cA)\in\stabdivS$ is a stability condition such that $\cO_S[1]\in\cA$, then $\cO_S[1]$ is stable with respect to $\sigma$.

\item
Let $S$ be a surface of Picard Rank $2$.
Assume that the cone of effective curves is generated by $C_1$ and $C_2$, that $C_1.G>0$ and that $C_1$ is the only irreducible curve in $S$ of negative self-intersection.
Then, $\cO_S[1]$ is only destabilized by $\cO_S(C_1)|_{C_1}$. 
\end{enumerate}
\end{pro}

\bibliographystyle{plainurl}

\bibliography{references}

\begin{thebibliography}{1}

\bibitem{ABL}
Daniele Arcara and Aaron Bertram.
\newblock Bridgeland-stable moduli spaces for {K}-trivial surfaces.
\newblock {\em JEMS}, 15(1):1--38, 2013.
\newblock (with an appendix by Max Lieblich).

\bibitem{ABCH}
Daniele Arcara, Aaron Bertram, Izzet Coskun, and Jack Huizenga.
\newblock The minimal model program for the {H}ilbert scheme of points on
  $\mathbb{P}^2$ and {B}ridgeland stability.
\newblock {\em Adv. Math.}, 235:580--626, 2013.

\bibitem{stabcondsonlocalp2}
Arend Bayer and Emanuele Macr\`{i}.
\newblock The space of stability conditions on the local projective plane.
\newblock {\em Duke Math. J.}, 160:263--322, 2011.

\bibitem{Beauville}
Arnaud Beauville.
\newblock {\em Complex Algebraic Surfaces (London Mathematical Society Student
  Texts)}, volume~34.
\newblock Cambridge University Press, 2nd edition, 1996.

\bibitem{biratgeomhilbsurfs}
Aaron Bertram and Izzet Coskun.
\newblock The birational geometry of the {H}ilbert scheme of points on
  surfaces.
\newblock In {\em Birational geometry, rational curves and arithmetic}, pages
  15--55. Springer, 2013.

\bibitem{stabcondsontricats}
Tom Bridgeland.
\newblock Stability conditions on triangulated cateogories.
\newblock {\em Ann. Math.}, 166:317--345, 2007.

\bibitem{MaciociaWalls}
Antony Maciocia.
\newblock Computing the walls associated to {B}ridgeland stability conditions
  on projective surfaces.
\newblock {\em Asian J. Math.}, 18(2), 2014.

\bibitem{Mart}
Cristian Martinez.
\newblock Duality, {B}ridgeland wall-crossing and flips of secant varieties.
\newblock URL: \url{http://arxiv.org/abs/1311.1183}.

\bibitem{stabcondsextrcontrs}
Yukinobu Toda.
\newblock Stability conditions and extremal contractions.
\newblock {\em Math. Ann.}, 357(2):631--685, 2013.

\end{thebibliography}

\end{document}